\documentclass[twoside,11pt,leqno]{amsart}
\usepackage{amssymb}
\swapnumbers
\usepackage[PostScript=dvips,nohug]{diagrams}
\newcounter{thmcounter}
\theoremstyle{plain}
\newtheorem{thm}{Theorem \Roman{thmcounter}}[section]
\newtheorem{fundthm}[thm]{The fundamental theorem of conformal geometry}
\newtheorem{prop}[thm]{Proposition}

\theoremstyle{definition}
\newtheorem{defn}[thm]{Definition}
\theoremstyle{remark}
\newtheorem{rem}[thm]{Remark}
\newtheorem{rems}[thm]{Remarks}
\newcommand{\proofof}[1]{\end{#1}\begin{proof}}
\newcommand{\emphdef}{\textit}
\numberwithin{equation}{section}

\newcommand{\acknowledge}{\medbreak\noindent\textbf{Acknowledgements. }}
\newcommand{\thismonth}{\ifcase\month\or
  January\or February\or March\or April\or May\or June\or
  July\or August\or September\or October\or November\or December\fi
  \space\number\year}
\newcommand{\bauth}[1]{\mbox{#1},}
\newcommand{\bart}[1]{\textit{#1},}
\newcommand{\bjourn}[3]{#1 \textbf{#2} (#3)}
\newcommand{\bbook}[1]{\textsl{#1},}

\newcommand{\bpp}[2]{pp.~#1--#2.}
\newcommand{\tcite}[1]{\textup{\cite{#1}}}
\DeclareMathAlphabet{\mathrmsl}{OT1}{cmr}{m}{sl}
\newcommand{\rssymb}[2]{\newcommand{#1}{{\mathrmsl{#2}}}}
\newcommand{\calsymb}[2]{\newcommand{#1}{{\mathcal{#2}}}}
\newcommand{\bbsymb}[2]{\newcommand{#1}{{\mathbb{#2}}}}
\newcommand{\lieoper}[2]{\newcommand{#1}{\mathop{\mathfrak{#2}}}}
\newcommand{\oper}[3][n]{\newcommand{#2}{\mathop{\mathrm{#3}}\ifx
  n#1\nolimits\else\limits\fi}}
\newcommand{\rsoper}[3][n]{\newcommand{#2}{\mathop{\mathrmsl{#3}}\ifx
  n#1\nolimits\else\limits\fi}}
\bbsymb\C{C}\bbsymb\F{F}\bbsymb\HQ{H}\bbsymb\N{N}\bbsymb\Q{Q}
\bbsymb\R{R}\bbsymb\T{T}\bbsymb\Z{Z}
\calsymb\cA{A}\calsymb\cB{B}\calsymb\cC{C}\calsymb\cD{D}\calsymb\cE{E}
\calsymb\cF{F}\calsymb\cG{G}\calsymb\cH{H}\calsymb\cI{I}\calsymb\cJ{J}
\calsymb\cK{K}\calsymb\cL{L}\calsymb\cM{M}\calsymb\cN{N}\calsymb\cO{O}
\calsymb\cP{P}\calsymb\cQ{Q}\calsymb\cR{R}\calsymb\cS{S}\calsymb\cT{T}
\calsymb\cU{U}\calsymb\cV{V}\calsymb\cW{W}\calsymb\cX{X}\calsymb\cY{Y}
\calsymb\cZ{Z}
\newcommand{\eps}{\varepsilon}

\oper\End{End}                      
\oper\Sym{Sym}\oper\Skew{Skew}\oper\Alt{Alt}
\oper\Weyl{Weyl}\oper\Cl{Cl}
\oper\Aut{Aut}                      
\oper\GL{GL}\oper\SL{SL}
\oper\CO{CO}\oper\On{O}\oper\SO{SO}
\oper\Pin{Pin}\oper\Spin{Spin}\oper\CSpin{CSpin}
\oper\Symp{Sp}\oper\Un{U}\oper\SU{SU}\oper\CU{CU}
\lieoper\der{der}                   
\lieoper\gl{gl}\lieoper\sgl{sl}
\lieoper\co{co}\lieoper\on{o}\lieoper\so{so}
\lieoper\symp{sp}\lieoper\un{u}\lieoper\su{su}\lieoper\csu{csu}
\newcommand{\abrack}[1]{[\mkern-3mu[#1]\mkern-3mu]}
\newcommand{\ip}[1]{\langle#1\rangle}

\newcommand{\lie}[1]{\mathfrak{#1}}

\newcommand{\restr}[1]{|_{\lower.5pt\hbox{${}_{#1}$}}}

\renewcommand{\leq}{\leqslant}
\newcommand{\dsum}{\oplus}               
\newcommand{\tsum}{\mathinner{\otimes}}  
\newcommand{\tens}{\mathbin{\otimes}}    
\newcommand{\setdif}{\smallsetminus}     %
\newcommand{\skwend}{\mathinner{\vartriangle}} 
\renewcommand{\wedge}{\mathinner{\mathchar"225E}}
\newcommand{\dual}{^{*\!}}
\newcommand{\Cinf}{\mathrm{C}^\infty}
\rsoper\Ad{Ad}
\rsoper\Vol{Vol}
\rsoper\kernel{ker}
\rsoper\image{im}
\rsoper\alt{alt}           
\rsoper\sym{sym}           
\rsoper\trace{tr}          
\rsoper\detm{det}          
\rsoper\Ricci{Ric}         
\rsoper\ricci{ric}         
\rsoper\scal{scal}         
\rsoper\divg{div}
\rssymb\Diff{Diff}
\rssymb\Vect{Vect}
\rssymb\ev{ev}
\rssymb\vol{vol}
\rssymb\ad{ad}
\rssymb\coad{coad}
\rssymb\iden{id}
\newcommand{\Hom}{\mathrm{Hom}}
\newcommand{\connect}{\#}

\newcommand{\low}{^{\vphantom x}}
\newcommand{\af}{{\mathit{af}}}
\newcommand{\sd}{{\raise1pt\hbox{$\scriptscriptstyle +$}}}
\newcommand{\asd}{{\raise1pt\hbox{$\scriptscriptstyle -$}}}
\newcommand{\sdasd}{{\raise1pt\hbox{$\scriptscriptstyle\pm$}}}
\newcommand{\conf}{\mathsf{c}}
\newcommand{\cip}{\ip}
\newcommand{\mult}{^{\scriptscriptstyle\times}}

\newcommand{\II}{\mathrm{I\kern-1.5ptI}}
\newcommand{\Hb}{\mathcal H}
\newcommand{\Vb}{\mathcal V}

\newcommand{\CP}[1]{\C P^{#1}}

\newcommand{\gmw}{\mathrmsl w}
\newcommand{\gmv}{\mathrmsl v}
\newcommand{\XH}{{\tilde X}}
\newcommand{\stB}{*\low_{\!B}}
\newcommand{\omB}{\omega\low_{\!B}}
\newcommand{\muB}{\mu\low_{\!B}}
\begin{document}
\title[Selfdual Einstein metrics and conformal submersions]
{Selfdual Einstein metrics\\ and conformal submersions}
\author{David M. J. Calderbank}
\address{Department of Mathematics and Statistics\\ University of Edinburgh\\
King's Buildings, Mayfield Road\\ Edinburgh EH9 3JZ\\ Scotland.}
\address{Centre de Math{\'e}matiques\\ Ecole Polytechnique\\ UMR 7640 du
  CNRS\\ F-91128 Palaiseau
Cedex\\ France.}
\email{davidmjc@maths.ed.ac.uk}
\date{\thismonth}
\begin{abstract}
  Weyl derivatives, Weyl-Lie derivatives and conformal submersions are
  defined, then used to generalize the Jones-Tod correspondence between
  selfdual $4$-manifolds with symmetry and Einstein-Weyl $3$-manifolds with an
  abelian monopole. In this generalization, the conformal symmetry is replaced
  by a particular kind of conformal submersion with one dimensional fibres.
  Special cases are studied in which the conformal submersion is holomorphic,
  affine, or projective. All scalar-flat K{\"a}hler metrics with such a
  holomorphic conformal submersion, and all four dimensional hypercomplex
  structures with a compatible Einstein metric, are obtained from solutions of
  the resulting ``affine monopole equations''.  The ``projective monopole
  equations'' encompass Hitchin's twistorial construction of selfdual Einstein
  metrics from three dimensional Einstein-Weyl spaces, and lead to an explicit
  formula for carrying out this construction directly. Examples include new
  selfdual Einstein metrics depending explicitly on an arbitrary holomorphic
  function of one variable or an arbitrary axially symmetric harmonic
  function.  The former generically have no continuous symmetries.
\end{abstract}
\maketitle
\section*{Introduction}

The aims of this paper are threefold: firstly, to advertise the notion of a
Weyl derivative both as a simple, but useful, tool in differential geometry,
and also as an object of study in its own right; secondly, to apply this tool
to the theory of conformal submersions, with particular attention to the case
of selfdual conformal $4$-manifolds; and thirdly to give explicit
constructions of selfdual Einstein metrics. The key discovery is a class of
conformal submersions with one dimensional fibres which admit a holomorphic
interpretation on the twistor space. This class includes the conformal
submersions generated by conformal vector fields and provides a natural
setting for a generalized Jones-Tod correspondence~\cite{JT}.

A Weyl derivative on a manifold $M$ is nothing more than a covariant
derivative on a real line bundle naturally associated to the differential
geometry of $M$. Weyl derivatives can be used to define Lie derivatives along
foliations with one dimensional leaves, generalizing the usual Lie derivative
along a vector field. They also occur naturally in the geometry of conformal
submersions. These two situations have in common the conformal submersions
with one dimensional fibres, to which most of this paper is devoted. I focus
in particular on the case that the total space is a selfdual $4$-manifold $M$
and define the notion of a selfdual conformal submersion. In Theorem
I~\eqref{JonesnTod}, the base $B$ of such a submersion is shown to be not just
conformal, but Einstein-Weyl, generalizing the Jones-Tod correspondence to a
context which also includes Hitchin's construction of selfdual Einstein
metrics with Einstein-Weyl conformal infinity~\cite{Hitchin3}.

The central part of the paper deals with special cases of this construction.
In Theorem II~\eqref{complex}, shear-free geodesic congruences on $B$ are
shown to correspond to antiselfdual complex structures on $M$ which are
invariant with respect to the conformal submersion (i.e., the submersion is
holomorphic), and the hypercomplex structures and scalar-flat K{\"a}hler metrics
arising in this way are identified (generalizing~\cite{CP2,GT}). In Theorem
III~\eqref{MT}, an antiselfdual complex structure is constructed from a
generic selfdual conformal submersion on a selfdual Einstein-Weyl
$4$-manifold---for a selfdual Einstein metric with a Killing field, this
reduces to the construction of Tod~\cite{Tod3}.

The notion of an affine conformal submersion is defined, characterized, and
shown, in Theorem IV~\eqref{affine}, to provide a method for constructing
selfdual spaces from coupled linear differential equations, which will be
called \emphdef{affine monopole equations}. In Theorem V~\eqref{SFK}, I show
that the generalized Jones-Tod constructions for scalar-flat K{\"a}hler and
hyperK{\"a}hler metrics arise in this way---this includes LeBrun's construction
of scalar flat metrics with Killing vector fields~\cite{LeBrun1} and the
construction of such metrics with homothetic vector fields~\cite{GT,CP2}. The
same ideas are applied to projective conformal submersions in Theorem
VI~\eqref{projective}, giving \emphdef{projective monopole equations}. This
time the differential equations are nonlinear: they are the $\SL(2,\R)$
Einstein-Weyl Bogomolny equations.

In the rest of the paper, I study conformal submersions on Einstein-Weyl
spaces, especially selfdual Einstein and locally hypercomplex $4$-manifolds.
In Theorem VII~\eqref{hitchin}, Hitchin's version of LeBrun's $\cH$-space
construction of selfdual Einstein metrics is characterized amongst selfdual
conformal submersions. In Theorem VIII~\eqref{biEW}, a larger class of
selfdual conformal submersions is obtained, from pairs of compatible selfdual
Einstein-Weyl structures. In Theorem IX~\eqref{einstein}, selfdual Einstein
metrics metrics with compatible hypercomplex structures are all found as
affine conformal submersions over hyperCR Einstein-Weyl spaces. This special
case of the Hitchin-LeBrun construction yields new selfdual Einstein metrics
with no continuous symmetries. Finally, in Theorem X~\eqref{final}, an
explicit formula for applying the Hitchin-LeBrun construction to any
Einstein-Weyl space is given.

The paper is organized as follows. In section~\ref{weyl}, Weyl derivatives and
Weyl-Lie derivatives are introduced. Although it is possible to present some
of the later results without this formalism, the proofs are simpler and more
natural when one makes systematic use of the affine space of Weyl derivatives
and the product rule for Weyl-Lie derivatives. Indeed, many of the results of
this paper would have been impossible to find (for the author at least)
without the geometric guidance provided by working in a gauge-independent
way. In order to familiarize the reader with this language, I have presented a
few simple applications of Weyl derivatives, and discussed carefully the
Weyl-Lie derivative on natural bundles. The key formula from this section is
the Weyl-Lie derivative of a torsion-free covariant derivative on a natural
bundle. In section~\ref{cg}, after recalling basic facts from conformal
geometry, I present another arrow in the Weyl geometer's quiver: the
linearized Koszul formula.

In the third section, the notion of a conformal submersion is defined, but
most of the local properties are studied within the more general framework of
conformal almost product structures. I prove a simple proposition which shows
that there is a canonically defined Weyl derivative in this setting, which
will be called the \emphdef{minimal Weyl derivative}. The main interest,
however, is in conformal submersions with one dimensional fibres, which may be
analyzed locally using the congruences and Weyl-Lie derivatives of
section~\ref{weyl}. I focus on this case for the second half of
section~\ref{cs} and characterize basic objects using the minimal Weyl-Lie
derivative.

The generalized Jones-Tod correspondence (Theorem I) is established in
section~\ref{JTC}.  In an earlier version of this paper~\cite{DMJC3}, my proof
followed closely the proof of the ``classical'' Jones-Tod correspondence given
in~\cite{CP1,CP2}.  However, Paul Gauduchon has recently obtained a cleaner
proof by exploiting more thoroughly the isomorphism between the bundle of
antiselfdual $2$-forms on the conformal $4$-manifold $M$ and the pullback of
the tangent bundle of the $3$-dimensional quotient $B$. This approach also has
the advantage of integrating nicely with the fact that invariant antiselfdual
complex structures on $M$ correspond to shear-free geodesic congruences on
$B$~\cite{CP2}, and so it is Gauduchon's proof that I follow here, adapted to
the context of conformal submersions.

This generalized Jones-Tod construction is difficult to apply, because the
generalized monopole equation and the defining equation for conformal
submersions are both nonlinear. This difficulty is partially overcome in
section~\ref{acs} by studying a special case in which the construction
linearizes: affine conformal submersions. In this case $M$ is an affine bundle
over $B$ such that the nonlinear connection and relative length scale induced
by the conformal structure are affine. A characterization of such submersions
is presented and the affine monopole equations are obtained. These linear
differential equations give constructions of scalar-flat K{\"a}hler and
hyperK{\"a}hler metrics extending work of LeBrun~\cite{LeBrun1} and Pedersen and
myself~\cite{CP2}. I illustrate this with some examples, taken from~\cite{CT}.
The following section extends these ideas to projective conformal submersions.
The projective monopole equations are not linear, but they will play a crucial
role in the final section.

Several of the results in this paper were motivated by twistor theory. I
explain this in section~\ref{twistor}, where I also prove that the
Hitchin-LeBrun construction, defined twistorially in~\cite{Hitchin3}, really
is a special case of the generalized Jones-Tod correspondence, and this
special case is characterized. The Hitchin-LeBrun construction, although
simple from a twistor point of view, has always been notoriously difficult to
carry out explicitly due to a lack of a direct construction: the known
examples are, to the best of my knowledge, those
of~\cite{LeBrun0,Hitchin3,HP2,CT}. The main result of the final portion of the
paper is Theorem X, which reduces the Hitchin-LeBrun construction to an
explicit formula for a selfdual Einstein metric in terms of an arbitrary
Einstein-Weyl structure. Firstly, though, a special case, Theorem IX, is
established using affine conformal submersions.  A key tool here is the
observation that two compatible Einstein-Weyl structures on a conformal
manifold define a conformal submersion, and in four dimensions, this
submersion is selfdual if the conformal structure is. This is proven in
section~\ref{EW}. Theorem IX then characterizes all selfdual Einstein
metrics admitting compatible hypercomplex structures. Explicit examples, with
no continuous symmetries, are given. They depend on an arbitrary holomorphic
function of one variable.

In the final section, a projective gauge transformation is applied to the
affine monopoles of Theorem IX. The result is a canonical solution of the
projective monopole equations which makes sense on any Einstein-Weyl space. In
Theorem X, I show that the explicit metric given by the induced projective
conformal submersion is Einstein. The resulting direct method for carrying out
the Hitchin-LeBrun construction is illustrated by a family of selfdual
Einstein metrics depending on an arbitrary axially symmetry harmonic function
on $\R^3$.

\acknowledge It is a pleasure to acknowledge the influence of Paul Gauduchon
on this work. Among numerous invaluable and fruitful conversations, I would
like to thank him in particular for explaining to me his approach to the
Jones-Tod correspondence and assisting with its adaptation to conformal
submersions. Discussions with him were also crucial for understanding and
characterizing the affineness of the nonlinear connection on a conformal
submersion. In addition, his suggestion that it would be worthwhile to give a
direct, rather than twistorial proof of Theorem III, led to a substantial
improvement in the presentation of the paper. I am grateful for the
hospitality of Ecole Polytechnique and a research grant from the EPSRC, which
made these discussions possible.

This paper has had a long evolution (for a primitive version see~\cite{DMJC3})
and conversations with many others have also been influential. I am grateful
to Florin Belgun and Andr{\'e} Moroianu for the idea to write the generalized
monopole equations as evolution equations. I would also especially like to
thank: Vestislav Apostolov, for discussing (bi-)hypercomplex structures and
Einstein metrics with me; Tammo Diemer, for pointing out, among other things,
the Koszul formula in Weyl geometry; Nigel Hitchin for suggesting that I
should look for affine structures compatible with conformal submersions; and
Henrik Pedersen, both for initiating my work on Weyl structures, conformal
submersions and selfdual Einstein metrics, and also for many discussions
relating to the Hitchin-LeBrun construction.

Snapshots from the affine space of Weyl derivatives in sections~\ref{EW}
and~\ref{SDE} were produced with the aid of Paul Taylor's Commutative Diagrams
package.

\section{Weyl derivatives and Weyl-Lie derivatives}\label{weyl}

If $V$ is a real $n$-dimensional vector space and $w$ any real number, then
the oriented one dimensional linear space $L^w=L^w(V)$ carrying the
representation $A\mapsto|\det A|^{w/n}$ of $\GL(V)$ is called the space
of \emphdef{densities of weight $w$} or \emphdef{$w$-densities}. It can be
constructed canonically as the space of maps
$\rho\colon(\Lambda^nV)\setdif0\to\R$ such that
$\rho(\lambda\omega)=|\lambda|^{-w/n}\rho(\omega)$ for all $\lambda\in\R\mult$
and $\omega\in(\Lambda^nV)\setdif0$.

The same construction can be carried out pointwise on any vector bundle $E$ to
give, for each $w\in\R$, the oriented real line bundle $L^w_E$ whose fibre at
$x$ is $L^w(E_x)$. Applying this to the tangent bundle gives the following
definition.
\begin{defn}
Suppose $M$ is any $n$-manifold. Then the \emphdef{density line bundle}
$L^w=L^w_{TM}$ of $M$ is defined to be the bundle whose fibre at $x\in M$
is $L^w(T_xM)$. Equivalently it is the associated bundle
$\GL(M)\times_{\GL(n)} L^w(n)$ where $\GL(M)$ is the frame bundle of $M$
and $L^w(n)$ is the space of $w$-densities of $\R^n$.
\end{defn}

The density bundles are oriented (hence trivializable) real line bundles, but
there is no preferred trivialization. Sections of $L=L^1$ may be thought of as
scalar fields with dimensions of length. This geometric dimensional analysis
may also be applied to tensors: the tensor bundle $L^w\tens(TM)^j\tens(T\dual
M)^k$ (and any subbundle, quotient bundle, element or section) is said to have
\emphdef{weight} $w+j-k$, or dimensions of [\textsl{length}]$^{w+j-k}$.  Note
that $L^{w_1}\tens L^{w_2}$ is canonically isomorphic to $L^{w_1+w_2}$ and
$L^0$ is the trivial bundle. When tensoring a vector bundle with some $L^w$
(or any real line bundle), I shall often omit the tensor product sign. Note
also that an orientation of $M$ may be viewed as a unit section of
$L^n\Lambda^nT\dual M$, defining an isomorphism between $L^{-n}$ and
$\Lambda^nT\dual M$.  A nonvanishing (usually positive) section of $L^1$ (or
$L^w$ for $w\neq0$) is called a \emphdef{length scale} or \emphdef{gauge} (of
weight $w$).

\begin{defn} A \emphdef{Weyl derivative} is a covariant derivative $D$ on
  $L^1$.  It induces covariant derivatives on $L^w$ for all $w$. The curvature
  of $D$ is a real $2$-form $F^D$ called the \emphdef{Faraday curvature}.
\end{defn}

If $F^D=0$ then $D$ is said to be \emphdef{closed}. There are then
local length scales $\mu$ with $D\mu=0$. If such a length scale exists
globally then $D$ is said to be \emphdef{exact}.  Conversely, a length scale
$\mu$ induces an exact Weyl derivative $D^\mu$ such that $D^\mu\mu=0$. Note
that $D^{c\mu}=D^\mu$ for any constant $c\neq0$. Weyl derivatives form an
affine space modelled on the linear space of $1$-forms, while closed and exact
Weyl derivatives are affine subspaces modelled on the linear spaces of closed
and exact $1$-forms respectively.

A \emphdef{gauge transformation} on $M$ is a positive function $e^f$ which
rescales a gauge $\mu\in\Cinf(M,L^w)$ to give $e^{wf}\mu$. It acts on Weyl
derivatives by $e^f\cdot D=e^f\circ D\circ e^{-f}=D-df$, so that $e^f\cdot
D^\mu=D^{e^f\mu}$ for $\mu\in\Cinf(M,L^1)$. If $D$ is any Weyl derivative,
then $D=D^\mu+\omega^\mu$ for the $1$-form $\omega^\mu=\mu^{-1}D\mu$, and
consequently, $\omega^{e^f\mu}=\omega^\mu+df$.

On an oriented manifold, Weyl derivatives may be viewed as a generalization of
volume forms, since the exact Weyl derivatives correspond to volume forms up
to constant multiples. For instance, the divergence of a vector field $X$ with
respect to a volume form $\vol$ is defined by $\cL_X\vol=(\divg X)\vol$. In
fact the divergence is naturally defined on vector field densities
$\Cinf(M,L^{-n}TM)$, and by twisting by a Weyl derivative $D$ on $L^n$ one can
define $\divg^DX$ for vector fields.

For another example, let $\Omega\in L^2\Lambda^2T\dual M$ be a weightless
$2$-form on $M^{2m}$ such that $\Omega^m$ is an orientation. Then $\Omega$
is nondegenerate and one would like to find a length scale $\mu$ such that
$\mu^{-2}\Omega$ is symplectic. If $2m>2$ this may not be possible.

\begin{prop}\textup{(cf.~\cite{Lee})} Let $M$ be an $n$-manifold
\textup($n=2m>2$\textup) and let $\Omega\in L^2\Lambda^2T\dual M$ be
nondegenerate. Then there is a unique Weyl derivative $D$ such that
$d^D\Omega$ is tracefree with respect to $\Omega$, in the sense that $\sum
d^D\Omega(e_i,e_i',.)=0$, where $e_i,e_i'$ are frames for $L^{-1}TM$ with
$\Omega(e_i,e_j')=\delta_{ij}$.
\proofof{prop} Pick any Weyl derivative $D^0$ and set $D=D^0+\gamma$ for some
$1$-form $\gamma$. Then $d^D\Omega=d^{D^0}\Omega+2\gamma\wedge\Omega$ and so
the traces differ by
\begin{equation*}\notag
2(\gamma\wedge\Omega)(e_i,e_i',.)
=2\gamma(e_i)\Omega(e_i',.)+2\gamma(e_i')\Omega(.,e_i)+
2\gamma\,\Omega(e_i,e_i')=2(n-2)\gamma.
\end{equation*}
Since $n>2$ it follows that there is a unique $\gamma$ such that
$d^D\Omega$ is tracefree.
\end{proof}

There are therefore two local obstructions to finding $\mu$ with
$d(\mu^{-2}\Omega)=0$, namely $d^D\Omega$ and $F^D$~\cite{Lee}. In four
dimensions $d^D\Omega$ automatically vanishes. In general if $d^D\Omega=0$
then $F^D\wedge\Omega=0$. In six or more dimensions this implies $F^D=0$, so
$D$ is closed, but it need not be exact. In four dimensions, however,
$D$ need not even be closed. This construction is of particular interest in
Hermitian geometry~\cite{Vaisman1}.

There is also a version of this in contact geometry. A \emphdef{contact
structure} on $M$ is codimension one subbundle $\Hb$ of $TM$ which is
maximally nonintegrable, in the sense that the Frobenius tensor
$\Omega\low_\Hb\colon\Lambda^2\Hb\to TM/\Hb$ is nondegenerate. In this
context, one defines a (slightly generalized) Weyl derivative to be a
covariant derivative on the real line bundle $TM/\Hb$. The following result is
then obtained.

\begin{prop} Let $M,\Hb$ be a contact manifold. Then there is a bijection
between complementary subspaces to $\Hb$ in $TM$ and Weyl derivatives on
$TM/\Hb$ such that the horizontal part of $F^D$ is tracefree with respect to
$\Omega\low_\Hb$.

\proofof{prop} Let $\eta\low_\Hb\colon TM\to TM/\Hb$ be the twisted
contact $1$-form whose kernel defines $\Hb$.  Then given any Weyl
derivative $D$, one can define $d^D\eta\low_\Hb$ and if $D=D^0+\gamma$
then $d^D\eta\low_\Hb=d^{D^0}\eta\low_\Hb+\gamma\wedge\eta\low_\Hb$.
Note that $d^D\eta\low_\Hb\restr\Hb$ is well defined, being equal to
$\Omega\low_\Hb$.  Since $\Omega\low_\Hb$ is nondegenerate, given $D$
there is a unique complementary subspace (spanned by $\xi$, say) such
that $d^D\eta\low_\Hb(\xi,.)=0$ and such a complementary subspace
fixes $D$ up to $1$-forms $\gamma$ with $\eta\low_\Hb\wedge\gamma=0$.

Now note that if $D=D^0+\mu^{-1}\eta\low_\Hb$ for a section
$\mu$ of $TM/\Hb$, then
$F^D\restr\Hb=F^{D^0}\restr\Hb+\mu^{-1}\Omega\low_\Hb$
and consequently, $D$ may be found uniquely with tracefree
horizontal Faraday curvature.
\end{proof}
This generalizes the fact that a contact form (which corresponds to a section
of $TM/\Hb$) defines a Reeb vector field complementary to the contact
distribution. Note also that it is perhaps more natural to work with
``horizontal'' covariant derivatives in this context, i.e., covariant
differentiation is only defined along directions in $\Hb$. Then the condition
on $F^D$ can be ignored, and the bijection becomes affine.

Much of the rest of the paper is concerned with the dual situation of one
dimensional subbundles of $TM$ (the complementary subspaces arising above
being an example). The integral manifolds of such a distribution define a
foliation of $M$ with one dimensional leaves.  This will be viewed as an
unparameterized version of a vector field by thinking of such a subbundle as
an inclusion $\xi\colon\Vb\to TM$ of a real line bundle $\Vb$, and hence as
a ``twisted'' vector field.  Sections of $\Vb$ correspond to vector fields
tangent to the foliation. It is therefore natural to consider covariant
derivatives on $\Vb$, which will again be referred to as Weyl derivatives.
I will also use the following terminology from relativity.

\begin{defn} A \emphdef{congruence} on a manifold $M$ is a nonvanishing
section $\xi$ of $\Vb^{-1}TM$ for some oriented real line bundle $\Vb$. It
defines an oriented one dimensional subbundle of $TM$ and hence a foliation
with oriented one dimensional leaves.
\end{defn}
No use will be made of the orientation of $\Vb$ in this section.

\begin{prop} Let $E$ be a vector bundle associated to the frame bundle
of $M$ with induced representation $\rho$ of $\lie{gl}(TM)$.  Then for a
congruence $\xi$, a Weyl derivative $D$ \textup(on $\Vb$\textup), and a
section $s$ of $E$, the formula
\begin{equation*}
\mu^{-1}\cL_{\mu\xi}s+\rho(\mu^{-1}D\mu\tens\xi)s
\end{equation*}
is independent of the choice of a nonvanishing section $\mu$
of $\Vb$, and will be called the \emphdef{Weyl-Lie derivative} $\cL^D_\xi s$
of $s$ along $\xi$.
\proofof{prop} This follows from the fact that for any vector field $X$
and function $f$, $\cL_{fX}s=f\cL_Xs-\rho(df\tens X)s$.
\end{proof}

Note that $\cL^D_\xi s$ is a section of $\Vb^{-1}E$. If $D$ is exact, then
trivializing $\Vb$ by a parallel section gives back the usual Lie derivative.
The dependence of the Weyl-Lie derivative on $D$ is clearly given by
$\cL^{D+\gamma}_\xi s=\cL^D_\xi s+\rho(\gamma\tens\xi)s$.  In particular, for
functions ($\rho$ trivial), $\cL^D_\xi f=df(\xi)$, which is a section of
$\Vb^{-1}$ independent of $D$.

In order to compute the Weyl-Lie derivative, it is convenient to find a
formula in terms of a torsion-free connection inducing covariant derivatives
$\nabla$ on any associated bundle $E$. In the case of the usual Lie
derivative, $\cL_KX=\nabla_KX-\nabla_XK$ for vector fields $X$,
and so $\cL_Ks=\nabla_Ks-\rho(\nabla K)s$ on sections of $E$. This
readily yields the following generalization to Weyl-Lie derivatives:
\begin{equation}\label{Lieform}
\cL^D_\xi s=\nabla_\xi s-\rho\bigl((D\tsum\nabla)\xi\bigr)s,
\end{equation}
where $(D\tsum\nabla)\xi$ denotes the \emphdef{twisted} or \emphdef{tensor sum}
covariant derivative of $\xi$ as a section of $\Vb^{-1}TM=\Vb^{-1}\tens TM$.
More generally $D\tsum\nabla$ will denote the twisted covariant derivative on
$\Vb^{-1}E$. Formula~\eqref{Lieform} immediately gives the following.

\begin{prop} Let $\xi$ be a congruence and $D$ be a Weyl derivative \textup(on
$\Vb$\textup). Then the Weyl-Lie derivative may be computed in terms of an
arbitrary torsion-free covariant derivative $\nabla$ as follows:
\begin{itemize}
\item For a $w$-density $\mu$, $\cL^D_\xi\mu=\nabla_\xi\mu-\frac
wn(\divg^{D\tsum\nabla}\xi)\mu$.
\item For a vector field $X$, $\cL^D_\xi X=\nabla_\xi
X-(D\tsum\nabla)\low_X\xi$.
\item For a $1$-form $\alpha$, $\cL^D_\xi\alpha=\nabla_\xi\alpha
+\alpha\bigl((D\tsum\nabla)\xi\bigr)=d\alpha(\xi,.)+D(\alpha(\xi))$.
\item For a $k$-form $\alpha$, $\cL^D_\xi\alpha=\iota\low_\xi
d\alpha+d^D(\iota\low_\xi\alpha)$.
\end{itemize}
\end{prop}
Here $\divg^{D\tsum\nabla}\xi$ denotes the trace of $(D\tsum\nabla)\xi$, i.e.,
the divergence has been twisted by $D$ on $\Vb^{-1}$ and $\nabla$ on $L^n$.

The above treatment only deals with the Weyl-Lie derivative for zero and first
order geometric objects (functions and sections of bundles associated to the
first order frame bundle). It may be extended to differential operators
(higher order geometric objects) using the product rule, but because the
Weyl-Lie derivative of a section of $E$ is not a section of $E$, the
differential operators have to be twisted by $D$. Consequently, some natural
differential operators may have nonzero Weyl-Lie derivative.  In particular
$(\cL^D_\xi d)\alpha=\cL^D_\xi(d\alpha)-d^D(\cL^D_\xi\alpha)=
F^D\wedge\iota\low_\xi\alpha$, i.e., the Weyl-Lie derivative
only commutes with exterior differentiation on functions in
general. Similarly, for the Lie bracket, $\cL_\xi^D[,](X,Y)=-F^D(X,Y)\xi$.  Of
course, this is compatible with the definition of $d$ in terms of $[,]$.

If $\nabla$ is any torsion-free covariant derivative then its Weyl-Lie
derivative $\cL^D_\xi\nabla=\cL^D_\xi\circ\nabla-(D\tsum\nabla)\circ\cL^D_\xi$
evaluates to
\begin{equation}
(\cL_\xi^D\nabla)\low_X=
R^\nabla_{\xi,X}+(D\tsum\nabla)\low_X\bigl((D\tsum\nabla)\xi\bigr)
\quad\in\:\Vb^{-1}\lie{gl}(TM).
\end{equation}
This acts on an associated bundle $E$ via the corresponding representation
$\rho$ of $\lie{gl}(TM)$. The vanishing of $\cL^D_\xi\nabla$ defines a notion
of invariance, along $\xi$, of $\nabla$ on $E$. Also note that
$\nabla_\xi-\cL^D_\xi=\rho\bigl((D\tsum\nabla)\xi\bigr)$ is a kind of Higgs
field.  If the Higgs field vanishes, I will say $\nabla$ is horizontal,
for reasons that will become clear later. In particular if $\nabla$ is
invariant and horizontal on $E$, then $\rho(R^\nabla_{\xi,X})=0$.

\section{Conformal geometry}\label{cg}

In the previous section, the term ``Weyl derivative'' was sometimes applied in
a generalized sense, when the oriented real line bundle was not necessarily
$L^1$. Such a distinction disappears when one introduces a conformal
structure.

\begin{defn} A \emphdef{conformal structure} $\conf$ on $M$ is a metric on
$L^{-1}TM$. Since the densities of $L^{-1}_xT_xM$ are canonically trivial,
it makes sense to require that this metric is \emphdef{normalized} in the
sense that it has determinant one.
\end{defn}

This is precisely what is needed in order to identify any oriented one
dimensional subbundle or quotient bundle of $TM$ with $L^1$.  In particular if
$M$ is conformal and $\xi$ is a congruence, then there is a unique oriented
isomorphism between $\Vb$ and $L^1$ such that $\xi\in\Cinf(M,L^{-1}TM)$ is a
weightless unit vector field. Thus a congruence will be viewed as an injective
linear map $\xi\colon L^1\to TM$ and $\Vb$ will denote its image.

A conformal structure may be viewed as a fibrewise inner product on $TM$ with
values in $L^2$: compatible Riemannian metrics therefore correspond
bijectively to length scales. As in the previous section, it is natural to
replace length scales with Weyl derivatives. Denoting the conformal inner
product of vector fields by $\cip{X,Y}\in\Cinf(M,L^2)$, the Koszul formula for
the Levi-Civita connection generalizes.

\begin{fundthm}\label{FTCG}\tcite{Weyl}
On a conformal manifold $M$ there is an affine bijection
between Weyl derivatives and torsion-free connections on $TM$
preserving the conformal structure. More explicitly, the torsion-free
connection on $TM$ is determined by the Koszul formula
\begin{equation*}\begin{split}
2\cip{D\low_XY,Z}&=
\quad D\low_X\,\cip{Y,Z}+D\low_Y\,\cip{X,Z}-D\low_Z\,\cip{X,Y}\\
&\quad+\cip{[X,Y],Z}-\cip{[X,Z],Y}-\cip{[Y,Z],X}.
\end{split}\end{equation*}
The corresponding linear map sends a $1$-form $\gamma$ to the
$\lie{co}(TM)$-valued $1$-form $\Gamma$ defined by $\Gamma_X=
\abrack{\gamma,X}=\gamma(X)\iden+\gamma\skwend X$, where $(\gamma\skwend
X)(Y)= \gamma(Y)X-\cip{X,Y}\gamma$.
\end{fundthm}
\begin{rem} Here, and elsewhere, I freely identify a $1$-form $\gamma$
with a vector field of weight $-1$ using the natural isomorphism $\sharp\colon
T\dual M\to L^{-2}TM$ given by the conformal structure. More generally, vector
fields of any weight are identified with $1$-forms of the same
weight. Similarly a skew linear map $J$ on $TM$ (of any weight) corresponds to
a $2$-form $\Omega\low_J$ (of the same weight) via
$J(X)=\sharp(\iota\low_X\Omega\low_J)$; this identifies $\gamma\skwend X$ with
$\gamma\wedge\cip{X,.}$. It will sometimes, but not always, be helpful to
maintain a notational distinction between a skew endomorphism and the
corresponding $2$-form. The bracket $\abrack{.\,,.}$ is part of an
algebraic Lie bracket on
$TM\dsum\lie{co}(TM)\dsum T\dual M$, and the same notation will be used for
the commutator bracket on $\co(TM)$.
\end{rem}

The curvature $R^D$ of $D$, as a $\lie{co}(TM)$-valued $2$-form, decomposes as
follows:
\begin{equation}\label{curv} R^D_{X,Y}=W\low_{X,Y}-\abrack{r^D(X),Y}+
\abrack{r^D(Y),X}.
\end{equation}
Here $W$ is the \emphdef{Weyl curvature} of the conformal structure, an
$\lie{so}(TM)$-valued $2$-form, and $r^D$ is a covector valued $1$-form,
the \emphdef{normalized Ricci tensor} of $D$.

If $\xi$ is a congruence, then the Weyl-Lie derivative of $D$ on $L^1$ reduces
to
\begin{equation} (\cL_\xi D)^{L^1}_X\mu
=F^D(\xi,X)\mu+\frac1n D^\xi_X(\divg^{D^\xi\tsum D}\xi)\mu,
\end{equation} where $D^\xi$ is the Weyl derivative used to define $\cL_\xi$.
Linearizing the Koszul formula gives a formula for the Weyl-Lie
derivative on other bundles:
\begin{equation*}\notag
(\cL_\xi D)\low_X=D\low_X(\cL_\xi\conf)+\alt\bigl(D(\cL_\xi\conf)(X)\bigr)
+\abrack{\cL_\xi D^{L^1},X}+\tfrac12\xi\skwend F^\xi(X)
+\tfrac12\cip{\xi,X}F^\xi,
\end{equation*}
where $F^\xi$ is the Faraday curvature of $D^\xi$: these $F^\xi$ terms come
from the Weyl-Lie derivative of the Lie bracket.
All the terms apart from the first belong the conformal Lie algebra
$\lie{co}(TM)$. This reflects the fact that if $\cL_\xi\conf=0$, the Weyl-Lie
derivative preserves conformal subbundles of natural vector bundles and so
extends to any bundle associated to the conformal frame bundle.

\section{Conformal submersions}\label{cs}

\begin{defn} Let $\pi\colon M\to B$ be a smooth surjective map between
conformal manifolds and let the \emph{horizontal bundle} $\Hb$ be the
orthogonal complement to the vertical bundle $\Vb$ of $\pi$ in $TM$.
Then $\pi$ will be called a \emphdef{conformal submersion} iff for all
$x\in M$, $d\pi_x\restr{\Hb_x}$ is a nonzero conformal linear map.
\end{defn}
It is not at all necessary to restrict attention to submersions. The base
could, for instance, be an orbifold, or be replaced altogether by the
horizontal geometry of a foliation.  However, since I am primarily interested
in the local geometry, I shall usually take the base to be a manifold.  A
bundle $\Hb$ complementary to $\Vb$ is often called a
\emphdef{\textup(nonlinear\textup) connection} on $\pi$: it is equivalently
determined by a projection $\eta\colon TM\to\Vb$, the \emphdef{connection
  $1$-form}.

\begin{prop}\label{confsub} If $\pi\colon M\to B$ is a submersion onto a
conformal manifold $B$, then conformal structures on $M$ making $\pi$ into a
conformal submersion correspond bijectively to triples $(\Hb,\conf^\Vb,\gmw)$,
where $\Hb$ is a connection on $\pi$, $\conf^\Vb$ is a conformal structure on
the fibres, and $\gmw\colon \pi^*L^1_{TB}\cong L^1_\Hb\to L^1_\Vb$ is a
\textup(positive\textup) isomorphism.
\proofof{prop} $d\pi\colon\Hb\to\pi^*TB$ is certainly an isomorphism, so
$L^2_\Hb\cong\pi^*L^2_{TB}$ and the conformal structure on $\Hb$ is obtained
by pullback. Combining this with $\conf^\Vb$ gives an $L^2_\Vb$ valued inner
product $\conf=\conf^\Vb\dsum\gmw^2\conf^\Hb$, which in turn determines an
isomorphism between $L_\Vb$ and $L_{TM}$ such that $\conf$ becomes a conformal
structure on $M$.
\end{proof}
The final ingredient $\gmw$ in this construction will be called a
\emphdef{relative length scale}, since it allows vertical and horizontal
lengths to be compared. The freedom to vary $\gmw$ generalizes the so called
``canonical variation'' of a Riemannian submersion, in which the fibre
metric is rescaled, while the base metric remains constant.

\medbreak
A natural generalization of conformal submersions, following Gray~\cite{Gray},
is a \emphdef{conformal almost product structure}. On a conformal manifold
$M$, this is a nontrivial orthogonal direct sum decomposition
$TM=\Vb\dsum_\perp\Hb$, i.e., $\Vb$ and $\Hb$ are nontrivial subbundles of
$TM$ and are orthogonal complements with respect to the conformal
structure. Although the roles of $\Vb$ and $\Hb$ are interchangeable, I will
call the corresponding tangent directions vertical and horizontal.

Given any Weyl derivative $D$, observe that that the vertical component
$(D\low_XY)^\Vb$ for $X,Y\in\Hb$ is tensorial in $Y$ and so defines a tensor
in $\Hb\dual\tens\Hb\dual\tens\Vb$. I will write
$(D\low_XY)^\Vb=\II_\Hb^D(X,Y)+\frac12\Omega\low_\Hb(X,Y)$, where $\II_\Hb^D$
is symmetric and $\Omega\low_\Hb$ is skew, and extend these fundamental forms
by zero to $T\dual M\tens T\dual M\tens TM$. Similarly the horizontal
component $(D\low_UV)^\Hb$ for $U,V\in\Vb$ defines tensors $\II_\Vb^D$ and
$\Omega\low_\Vb$ in $\Vb\dual\tens\Vb\dual\tens\Hb\leq T\dual M\tens T\dual
M\tens TM$. Since $D$ is torsion-free, $\Omega\low_\Hb$ and $\Omega\low_\Vb$
are the Frobenius tensors of the distributions $\Hb$ and $\Vb$, which vanish
iff the distributions are tangent to foliations. On the other hand, the
fundamental forms $\II_\Hb^D$ and $\II_\Vb^D$ do depend on $D$ and can be used
to find a distinguished Weyl derivative $D^0$.

\begin{prop} Suppose $M$ is conformal with a conformal almost product
structure $TM=\Vb\dsum_\perp\Hb$. Then there is a unique Weyl derivative $D^0$
such that $\Vb$ and $\Hb$ are minimal, in the sense that the fundamental
forms, denoted $\II_\Hb^0$ and $\II_\Vb^0$, are tracefree. It may be
computed from an arbitrary Weyl derivative $D$ via the formula:
\begin{equation*}
D^0=D-\frac{\trace\II_\Hb^D}{\dim\Hb}-\frac{\trace\II_\Vb^D}{\dim\Vb}.
\end{equation*}
\proofof{prop} Observe that if $\tilde D=D+\gamma$ then
$\cip{\trace\II_\Hb^{\tilde D},U}= \cip{\trace
\II_\Hb^D,U}+(\dim\Hb)\gamma(U)$ for all $U\in\Vb$, and similarly for
$\II_\Vb$. This shows that the formula for $D^0$ is independent of $D$.
Substituting $D=D^0$ shows that $\II_\Hb^0$ and $\II_\Vb^0$ are tracefree.
\end{proof}
I shall refer to $D^0$ as the \emphdef{minimal Weyl derivative} of $\xi$; this
usage is consonant both with minimal submanifolds and minimal coupling. The
minimal Weyl derivative need not be closed: I shall write $F^0$ for its
Faraday curvature. Also I denote the minimal Weyl-Lie derivative along $\xi$
by $\cL^0_\xi$.

\begin{rem}
The curvature of $D^0$ on $TM$ can be related to the curvatures of horizontal
and vertical connections on $\Hb$ and $\Vb$. One defines a horizontal
connection on $\Hb$ by $D^\Hb_XY=(D^0_XY)^\Hb$ for $X,Y\in\Hb$; similarly
$D^\Vb_UV=(D^0_UV)^\Vb$ for $U,V\in\Vb$. The (modified) curvature of $D^\Hb$
is defined by
\begin{equation*}
R^\Hb_{X,Y}Z=D^\Hb_XD^\Hb_YZ-D^\Hb_YD^\Hb_XZ-D^\Hb_{[X,Y]^\Hb}Z
-\bigl[[X,Y]^\Vb,Z\bigr]{}^\Hb,
\end{equation*} where $X,Y,Z\in\Hb$. The definition of $R^\Vb$
is analogous, and O'Neill-type formulae~\cite{ONeill} relating $R^{D^0}$ to
$R^\Hb$ and $R^\Vb$ follow directly as in~\cite{Gray}.
\end{rem}

Some of the properties of the minimal Weyl derivative may be elucidated by
comparing it with partial Lie derivatives. If $X$ is a horizontal vector field
and $U$ is vertical, then $[U,X]^\Hb$ is tensorial in $U$, which defines a
partial covariant derivative on $\Hb$ in vertical directions.  This extends
naturally to horizontal forms in $\Lambda^k\Hb\dual$ and to densities in
$L^1_\Hb$. One says a horizontal vector field, form, or density is invariant
if its partial covariant derivative along $\Vb$ vanishes; if $\Vb$ is a
tangent to the fibres of a submersion, this means the horizontal vector field,
form, or density is basic.

Similarly, there is a partial covariant derivative on $\Vb$,
$\Lambda^k\Vb\dual$ and $L^1_\Vb$ in horizontal directions. Since the
conformal structure identifies $L^1_\Hb$ with $L^1_\Vb$, putting these
together gives a Weyl derivative.  In order to verify that this is the minimal
Weyl derivative defined above, introduce an arbitrary Weyl connection $D$ so
that $[U,X]^\Hb=(D\low_UX)^\Hb+(D\low_UX)^\Vb$ and the second term is
tensorial: since $\cip{D\low_UX,V}=-\cip{X,D\low_UV}$ this tensor is
essentially $\II_\Vb^D+\frac12\Omega\low_\Vb$. Therefore the vertical partial
connections on $L^1_\Hb$ induced by $[U,X]^\Hb$ and $(D\low_UX)^\Hb$ agree if
and only if $\trace\II_\Vb^D=0$.  Similarly the horizontal partial connections
on $L^1_\Vb$ induced by $[X,U]^\Vb$ and $(D\low_XU)^\Vb$ agree if and only if
$\trace\II_\Hb^D=0$.

According to this discussion, the following definition for densities is
compatible with the identification of $L^1$ with $L^1_\Hb$ and $L^1_\Vb$.

\begin{defn} A density $\mu\in\Cinf(M,L^1)$ on a conformal manifold $M$ with
a conformal almost product structure $(\Vb,\Hb)$ is \emphdef{invariant along
$\Vb$} iff $D^0_U\mu=0$ for all vertical $U$, and \emphdef{invariant along
$\Hb$} iff $D^0_X\mu=0$ for all horizontal $X$.
\end{defn}

I specialize now to the case that $\Vb$ is one dimensional and oriented. Then
the positively oriented weightless unit vector field $\xi$ spanning
$L^{-1}\Vb\leq L^{-1}TM$ is a congruence and the minimal Weyl derivative $D^0$
is characterized by $D^0_\xi\xi=0$ and $\trace D^0\xi=0$ (note that
$\cip{D^0\xi,\xi}=0$ since $\xi$ has unit length). The formula for computing
$D^0$ reduces to $D^0=D-\frac1{n-1}(\divg^D\xi)\xi+(d^D\xi)(\xi,.)$. Also
note that $\Omega\low_\Vb$ and $\II_\Vb^0$ both vanish, so $\Omega$ and
$\II^0$ will denote the fundamental forms of $\Hb$.

\begin{prop} Let $\xi$ be a congruence with minimal Weyl derivative $D^0$
and let $D$ be an arbitrary Weyl derivative.
\begin{enumerate}
\item $\frac12\cL^0_\xi\conf=\sym_0(D^0\tsum D)\xi=\sym D^0\xi=\II^0$.
\item $D^0$ is exact iff $\xi=K/|K|$ for some \textup(nonvanishing\textup)
vector field $K$ which is divergence-free with respect to the metric
$g=|K|^{-2}\conf$.
\item If $D^0$ is exact and $\cL^0_\xi\conf=0$ then $K$ is a conformal
vector field, and hence is a Killing field of the metric $g=|K|^{-2}\conf$.
\end{enumerate}
Conversely if $K$ is a nonvanishing conformal vector field then
$\xi=K/|K|$ is a congruence with $\cL^0_\xi\conf=0$ and $D^0|K|=0$.
\proofof{prop} For the first part, note that for any vector fields $X,Y$,
\begin{align*}
(\cL^0_\xi\conf)(X,Y)&=\cL^0_\xi\cip{X,Y}-\cip{\cL^0_\xi X,Y}
-\cip{X,\cL^0_\xi Y}\\
&=D\low_{\smash[b]\xi}\cip{X,Y}-\tfrac2n(\divg^{D^0\tsum
D}\xi)\cip{X,Y}\\
&\quad-\cip{D\low_{\smash[b]\xi} X-(D^0\tsum D)\low_X\xi,Y}
-\cip{X,D\low_{\smash[b]\xi}Y-(D^0\tsum D)\low_Y\xi}\\
&=\cip{(D^0\tsum D)\low_X\xi,Y}+\cip{X,(D^0\tsum D)\low_Y\xi}
-\tfrac2n(\divg^{D^0\tsum D}\xi)\cip{X,Y}.
\end{align*}
This is $2(\sym_0 (D^0\tsum D)\xi)(X,Y)=2(\sym D^0\xi)(X,Y)$, since $D$ was
arbitrary. If either $X$ or $Y$ is parallel to $\xi$, $(\sym D^0\xi)(X,Y)$
vanishes automatically, because $D^0_\xi\xi=0=\cip{D^0\xi,\xi}$. On the other
hand, if $X$ and $Y$ are orthogonal to $\xi$, it is equal to
$\frac12\cip{\xi,D^0_XY+D^0_YX}= \cip{\xi,\II^0(X,Y)}$.

For the second and third parts, observe that an exact $D^0$ preserves a length
scale $\mu$. Then $K=\mu\xi$ and $D^0$ is the Levi-Civita connection of
$g=\mu^{-2}\conf$.
\end{proof}

I shall now assume that the congruence $\xi$ is tangent to the one dimensional
fibres of a submersion over a manifold $B$; this is always true locally. A
horizontal vector field, form or density is then basic if it is invariant in
the sense above. The Weyl-Lie derivative $\cL^0_\xi$ provides an efficient way
to characterize such basic objects.

\begin{prop} Let $\xi$ be a congruence with minimal Weyl derivative
$D^0$ generating a submersion of $M$ over $B$. Then a horizontal vector field
$X$ is basic iff $\cL_\xi^0X=0$. Similarly a horizontal form $\alpha$ is basic
iff $\cL_\xi^0\alpha=0$. Finally a density $\mu$ is basic iff
$\cL^0_\xi\mu=0$.

\proofof{prop} It suffices to check that $\cL_\xi^0X=0$ is equivalent to
$[U,X]$ being vertical for all vertical $U$. If $U=\lambda\xi$ then away from
the zero set of $\lambda$, $D^\lambda=D^0+\gamma$ for some $1$-form $\gamma$
and so $[U,X]=\cL_{\lambda\xi}X=\lambda\cL^{D^\lambda}_\xi X= \lambda\cL^0_\xi
X+\gamma(X)\lambda\xi$. Hence $\cL^0_\xi X$ is vertical if and only if $[U,X]$
is vertical for all vertical $U$. For the vertical component observe that
$\cip{\xi,\cL^0_\xi X}=-\cL^0_\xi\conf(\xi,X)$ which vanishes by the previous
Proposition. The result for forms follows from the product rule.

For densities, $\cL^0_\xi\mu=D^0_\xi\mu$, and this means $\mu$ is basic as a
section of $L^1_\Hb$.
\end{proof}

The product rule means that other basic objects are characterized by vanishing
Weyl-Lie derivative. For instance the submersion is conformal (i.e., the
horizontal conformal structure is basic) iff $\cL^0_\xi\conf=0$.

Similarly a connection $\nabla$ on $E$ (associated to the frame bundle, or the
conformal frame bundle if $\cL^0_\xi\conf=0$) is basic if it is horizontal and
invariant with respect to $\cL^0_\xi$. In fact invariance suffices for the
horizontal part of a connection $\nabla$ to descend to $B$, but the pullback
connection is then
$\nabla-\cip{\xi,.}\tens\rho\bigl((D^0\tsum\nabla)\xi\bigr)$.

The horizontal part of a Weyl derivative $D$ on $M$ is basic iff
$0=\cL^0_{\smash[b]\xi} D=F^D(\xi,.)+\frac1n D^0(\divg^{D^0\tsum D}\xi)$. If
$D$ is horizontal this reduces to $F^D(\xi,.)$.  In particular $D^0$ itself is
basic iff $F^0(\xi,.)=0$.  If $\cL_\xi^0\conf=0$, the only nonzero fundamental
form is $\Omega=\Omega\low_\Hb$:
\begin{equation*}
\cip{\Omega(X,Y),\xi}=2\cip{D^0_XY,\xi}=-2\cip{D^0_X\xi,Y}=
-(d^0\xi)(X,Y).
\end{equation*}
Note that $\cL^0_\xi(d^0\xi)=\iota\low_\xi(d^0)^2\xi =F^0-\xi\wedge
\bigl(F^0(\xi,.)\bigr)$ and so $\Omega$ is basic iff $\xi\wedge F^0=0$.

Note also that the linearized Koszul formula with respect to the minimal Weyl
derivative of a conformal submersion $\xi$ reduces to:
\begin{equation}\label{linKoz}
(\cL^0_\xi D)\low_X=
\abrack{\cL^0_\xi D^{L^1},X}+\tfrac12\xi\skwend F^0(X)
+\tfrac12\cip{\xi,X}F^0.
\end{equation}
Hence invariance on $L^1$ does not imply invariance on other natural bundles.

\section{The Jones-Tod correspondence}\label{JTC}

Suppose that $\xi$ is a congruence on an oriented conformal $4$-manifold $M$
defining a conformal submersion $\pi$ onto a manifold $B$. Let $D^0$ be the
minimal Weyl derivative of $\xi$, define $\omega=-(*{d^D\xi})(\xi,.)$ (which
can be computed using any Weyl derivative $D$) and let
$D^{sd}=D^0+\frac12\omega$ and $D^B=D^0+\omega$.

\begin{rems} The definition of $\omega$ uses the natural extension of the
star operator to $2$-forms of any weight. The star operator on $1$-forms (of
any weight) is a $3$-form of the same weight defined by
$\iota\low_X{*\alpha}=*(\cip{X,.}\wedge\alpha)$ for any vector field $X$.
In general the star operator on any manifold will be defined in terms of the
orientation $*1$ so that a similar relation holds between the star
operator, wedge product and interior multiplication, with no signs. As
remarked in~\cite{CP2}, this is more convenient in computations than the
usual choice. Note that $*^2=+1$ in four dimensions, whereas $*_{\!B}^2=-1$
in three dimensions.

Since the star operator is an involution on $2$-forms in four dimensions,
$\Lambda^2T\dual M=\Lambda^2_+T\dual M\dsum\Lambda^2_-T\dual M$.  The selfdual
and antiselfdual parts of a $2$-form are denoted $F=F_++F_-$.  A skew
endomorphism $J$ of $TM$ may be identified with a weightless $2$-form
$\Omega\low_J\in L^2\Lambda^2T\dual M\cong L^{-2}\Lambda^2TM$ via
$\Omega\low_J(X,Y)=\cip{JX,Y}$ and $J$ is said to be \emphdef{selfdual} or
\emphdef{antiselfdual} if $\Omega\low_J$ is. If $J$ is antiselfdual, then
$\Omega\low_J=\eta\wedge J\eta-*\eta\wedge J\eta$ for any weightless unit
$1$-form $\eta$. It follows that, for a $1$-form $\alpha$ (of any weight),
${*\alpha}=J\alpha\wedge\Omega\low_J$ for any antiselfdual endomorphism $J$
with $J^2=-\iden$. Note that if $F$ is a selfdual $2$-form, then $X\wedge
{\,F(Y)}-Y\wedge {\,F(X)}$ is also selfdual.
\end{rems}

\begin{prop}\label{selfdual} $\cip{(D^0\tsum D^{sd})\xi,.}$  is a selfdual
$2$-form of weight $-1$.
\proofof{prop} In terms of an arbitrary Weyl derivative $D$,
\begin{align*}
D^0&=D-\tfrac14(\divg^{D^0\tsum D}\xi)\xi+\tfrac12(d^{D^0\tsum D}\xi)(\xi,.)\\
\omega&=-({*d^{D^0\tsum D}\xi})(\xi,.)\\
\tag*{and so}
D^{sd}&=D-\tfrac14(\divg^{D^0\tsum D}\xi)\xi
+\tfrac12(d^{D^0\tsum D}\xi)(\xi,.)-\tfrac12({*d^{D^0\tsum D}\xi})(\xi,.).
\end{align*}
Substituting $D=D^{sd}$ into this formula gives the result.
\end{proof}
\begin{rem} This property clearly characterizes $D^{sd}$. One can
characterize $D^B$ in a similar way by the vanishing of the trace of $D^B\xi$
and the selfduality of $\alt D^B\xi$.
\end{rem}

\begin{prop}\label{sdbasic} The Weyl derivative $D^B$ is basic iff $D^0$ has
selfdual Faraday curvature.
\proofof{prop} Since $\omega(\xi)=0$, $D^B$ is basic iff
$F^B(\xi,.)=0$, where $F^B$ is the Faraday curvature of $D^B$.
Now $D^B=D^0+\omega$ and so 
\begin{equation*}
F^B(\xi,.)=F^0(\xi,.)+d\omega(\xi,.)=F^0(\xi,.)+\cL^0_\xi\omega.
\end{equation*}
Writing $\omega=-(*d^0\xi)(\xi,.)=-{*(\xi\wedge d^0\xi)}$ yields
\begin{equation*}
\cL^0_\xi\omega=-{*(\xi\wedge\cL^0_\xi d^0\xi)}=-{*(\xi\wedge F^0)}
=-(*F^0)(\xi,.).
\end{equation*}
So $F^B(\xi,.)=(F^0-{*F^0})(\xi,.)$. Since $F^0-{*F^0}$ is antiselfdual,
this contraction with $\xi$ vanishes iff $F^0={*F^0}$.
\end{proof}
When $D^0$ has selfdual Faraday curvature, the conformal submersion is said to
be \emphdef{selfdual}. In this case $D^B$ is a basic Weyl derivative on
$L^1\cong\pi^*L^1_B$ and the induced Weyl structure on $B$ is sometimes called
the \emphdef{Jones-Tod Weyl structure}. It follows from the Koszul formula
that the induced Weyl connection on $TB$ pulls back to the conformal
connection on $\Hb\cong\pi^*TB$ given by the horizontal part of the Weyl
connection induced by $D^B$ on $TM$.  The same observation holds for
$L^{-1}\Hb\cong\pi^*L^{-1}_BTB$.

Now observe that the map $\Omega\low_J\mapsto J\xi$ is an isomorphism from
$L^2\Lambda^2_-T\dual M$ to $L^{-1}\Hb$ with inverse
$\chi\mapsto\xi\wedge\chi-{*(\xi\wedge\chi)}\in L^{-2}\Lambda^2_-TM\cong
L^2\Lambda^2_-T\dual M$. If $J^2=-\iden$ then $|J\xi|=1$
and so this isomorphism is an isometry up to a constant multiple
(conventionally, $|\Omega\low_J|^2=2$ when $J^2=-\iden$). Also note that
$\abrack{J_1,J_2}\xi=-2{*(\xi\wedge J_1\xi\wedge J_2\xi)}$.

The action of $\lie{co}(TM)$ on antiselfdual endomorphisms is by commutator,
and so only the antiselfdual part contributes, since selfdual and antiselfdual
endomorphisms commute. Therefore Proposition~\ref{selfdual} shows that
$\cL^0_\xi J=D^{sd}_\xi J$, i.e., $D^{sd}$ is horizontal on
$L^2\Lambda^2_-T\dual M$.  The linearized Koszul formula may be used to show
that if $\xi$ is selfdual, then $D^{sd}$ is invariant on $L^2\Lambda^2_-T\dual
M$. More precisely,
\begin{equation*}
\cL^0_\xi D^{sd}_X=\tfrac12 F^0(\xi,X)\iden
+\tfrac12\xi\skwend F^0(X)-\tfrac12X\skwend F^0(\xi)+\tfrac12\cip{\xi,X}F^0,
\end{equation*}
which has selfdual skew part. Hence $D^{sd}$ is basic on
$L^2\Lambda^2_-T\dual M\cong\pi^*L^{-1}_BTB$. The following Proposition
identifies it with $D^B$ (which gives another proof of invariance).

\begin{prop}\label{identification}
If $J$ is an antiselfdual endomorphism, then for any horizontal
vector field $X$, $\bigl(D^B_X(J\xi)\bigr){}^\Hb=(D^{sd}_X J)\xi$.
\proofof{prop} Write $\chi=J\xi$ so that
$\Omega\low_J=\xi\wedge\chi-{*(\xi\wedge\chi)}$. Then
\begin{align*}
D^{sd}_X\Omega\low_J&=D^{sd}_X\xi\wedge\chi+\xi\wedge D^{sd}_X\chi
-*(D^{sd}_X\xi\wedge\chi+\xi\wedge D^{sd}_X\chi).\\
\intertext{Now $D^{sd}=D^0+\frac12\omega=D^B-\frac12\omega$ and so}
D^{sd}_X\xi&=\tfrac12\,{*(X\wedge\xi\wedge\omega)},\\
D^{sd}_X\chi&=D^B_X\chi+\tfrac12\omega(\chi)X-\tfrac12\cip{\chi,X}\omega.\\
\tag*{Therefore}
D^{sd}_X\xi\wedge\chi&=-\tfrac12\cip{\chi,X}\,{*(\xi\wedge\omega)}
+\tfrac12\omega(\chi)\,{*(\xi\wedge X)}\\
\tag*{and}
\xi\wedge D^{sd}_X\chi&=\xi\wedge D^B_X\chi+\tfrac12\omega(\chi)\xi\wedge X
-\tfrac12\cip{\chi,X}\xi\wedge\omega
\end{align*}
which gives $\xi\wedge D^{sd}_X\chi-*(D^{sd}_X\xi\wedge\chi)=\xi\wedge
D^B_X\chi$.

Taking the antiselfdual part and contracting with $\xi$ completes the proof.
\end{proof}

A generalized Jones-Tod correspondence follows readily from these
observations, following an approach due to Gauduchon. Recall~\cite{Hitchin3}
that a Weyl connection is said to be \emphdef{Einstein-Weyl} iff the symmetric
traceless part of its Ricci tensor vanishes.  \addtocounter{thmcounter}{1}
\begin{thm}\label{JonesnTod} Suppose $(M,\conf)$ is an oriented conformal
$4$-manifold and $\xi$ is a selfdual conformal submersion over a manifold $B$.
Then $D^B=D^0+\omega$ is Einstein-Weyl on $B$ if and only if $\conf$ is
selfdual.  \proofof{thm} Since $\xi$ is selfdual, $D^B$ descends to a Weyl
connection on $B$. If $\pi^*D^B$ denotes the pullback of $D^B$ to
$\pi^*L^{-1}_BTB\cong L^{-1}\Hb$ then Proposition~\ref{identification} implies
that $D^{sd}J=(\pi^*D^B)(J\xi)$ for any antiselfdual endomorphism $J$, and so
$\abrack{R^{sd}_{X,Y},J}\xi=(\pi^*R^B)\low_{X,Y}(J\xi)$, where $R^B$ is the
curvature of $D^B$ on $B$ and $X,Y$ are arbitrary vector fields. Here I have
used the fact that $D^{sd}$ is horizontal, and so $D^{sd}_\xi J=\cL^0_\xi J
=(\pi^*D^B)\low_\xi(J\xi)$ by definition of pullback; horizontality and the
definition of pullback likewise imply that $\abrack{R^{sd}_{X,Y},J}$ and
$(\pi^*R^B)\low_{X,Y}(J\xi)$ vanish if $X$ or $Y$ is vertical. (Recall that
$\cL^0_\xi D^{sd}=R^{sd}_{\xi,X}+(D^0\tsum D^{sd})\low_X(D^0\tsum
D^{sd})\xi$.)

Let $R^{sd,\asd}\colon\Lambda^2TM\to L^2\Lambda^2T\dual_-M$ denote the
antiselfdual part of $R^{sd}$ and let $R^{B,0}\colon\Lambda^2TB\to
L^2_B\Lambda^2T\dual B$ denote the skew part of $R^B$. Then, omitting
pullbacks, $R^{sd,\asd}(J)\xi=-\frac12{\stB R^{B,0}({\stB J\xi})}$, since
$R^{sd,\asd}(\xi\wedge X)=0$. The symmetric traceless part of $J\mapsto
R^{sd,\asd}(J)$ is the antiselfdual Weyl tensor $W^\asd$ and the symmetric
traceless part of $\chi\mapsto{\stB R^{B,0}({\stB \chi})}$ is the
symmetric traceless Ricci tensor of $D^B$, which proves the theorem.
\end{proof}
An explicit formula for the relationship between $r^{sd}$
and $r^B$ will be useful later.
\begin{prop}\label{ricci} Suppose that $\xi$ is a selfdual conformal
submersion from a selfdual space $M$ to an Einstein-Weyl space $B$. Then
$F^{sd}_\asd=\frac14(F^B-{*F^B})$ and
\begin{equation*}
\sym r^{sd}=\tfrac1{12}\scal^B(\iden-2\xi\tens\xi)
+\tfrac14\bigl({\stB F^B}\tens\xi+\xi\tens{\stB F^B}\bigr).
\end{equation*}
\proofof{prop} Let $X,Y$ be basic vector fields. Since
$\cip{R^{sd,\asd}(\xi\wedge X),\xi\wedge Y}=0$, it follows that
\begin{equation*}
r^{sd}(X,Y)+r^{sd}(\xi,\xi)\cip{X,Y}
+{*\bigl(\xi\wedge r^{sd}(\xi,.)\wedge X\wedge Y\bigr)}=0.
\end{equation*}
On the other hand, since
$\cip{R^{sd}({*\xi\wedge X}),\xi\wedge Y-{*(\xi\wedge Y)}}
=\cip{R^{B,0}(\stB JX),\stB  JY}$, it follows that
\begin{multline*}
r^{sd}(X,Y)-(\trace\low_{\Hb}r^{sd})\cip{X,Y}
-{*\bigl(\xi\wedge r^{sd}(.,\xi)\wedge X\wedge Y\bigr)}\\
=r^B(X,Y)-(\trace r^B)\cip{X,Y}
=-\tfrac12F^B(X,Y)-\tfrac16\scal^B\cip{X,Y}.
\end{multline*}
The stated formulae follow easily from these.
\end{proof}
Note that $F^{sd}=\frac12(F^0+F^B)$ and so the formula for $F^{sd}_\asd$
follows immediately from the fact that $F^0_\asd=0$. On the other hand
$F^0$ and $F^{sd}_\sd$ are not basic in general.

The form of the Jones-Tod construction stated in Theorem I gives a procedure
for constructing Einstein-Weyl spaces from selfdual spaces. For the inverse
construction, the following reformulation is useful.
\begin{prop}\label{char} Suppose that $(M,\conf)$ is an oriented conformal
$4$-manifold, that $\xi$ is a conformal submersion over an Einstein-Weyl
manifold $B$, and that $D^0=\pi^*D^B-\omega$, where $D^B$ is the Weyl
derivative on $L^1_B$ and $\omega=-(*{d^D\xi})(\xi,.)$
\textup(computed using any Weyl derivative $D$\textup). Then $(M,\conf)$ is
selfdual and $\xi$ is selfdual.
\end{prop}
This follows immediately from Theorem I and Proposition~\ref{sdbasic}:
$\pi^*D^B-\omega$ has selfdual Faraday curvature since $\pi^*D^B$ is basic.
Therefore, an inverse Jones-Tod construction will be obtained if the equation
$D^0=\pi^*D^B-\omega$ can be interpreted as an equation on $B$.  In order to
do this, recall from Proposition~\ref{confsub}, that a conformal structure on
a fibre bundle $\pi\colon M\to B$ over a conformal manifold $B$ is determined
by a connection $1$-form $\eta\colon TM\to\Vb$ and a relative length scale
$\gmw\colon\pi^*L^{1}_B\to\Vb$, where I have assumed the fibres are
oriented and one dimensional so that $L^1_\Vb=\Vb$. Choosing a fibre
coordinate $t$ identifies $M$ locally with $B\times\R$, providing a
trivialization of $\Vb$ and a flat connection $1$-form $dt$. In these terms
$\eta=dt+A$ for $A\in\Cinf(B\times\R,\pi^*T\dual B)$ and
$\gmw\in\Cinf(B\times\R,\pi^*L^{-1}_B)$.

The inverse Jones-Tod construction can now be formulated as a nonlinear
evolution equation on $B$.
\begin{prop}\label{evolve}  Let $(M,\conf)$ be a selfdual conformal
$4$-manifold with a selfdual conformal submersion $\xi$ over an Einstein-Weyl
space $(B,\conf\low_B,D^B)$. Then $M\to B$ is locally conformal to $\pi\colon
B\times\R\to B$ with conformal structure
\begin{align*}
\conf=\pi^*\conf\low_B &+ \gmw^{-2}(dt+A)^2\\
\tag*{\textit{where}} \stB (D^B\gmw+\dot A\gmw-A\dot\gmw)&=dA+\dot A\wedge A,
\end{align*}
for $\gmw\in\Cinf(B\times\R,L^{-1}_B)$ and $A\in\Cinf(B\times\R,T\dual B)$.
Here $t$ is the fibre coordinate on $B\times\R$ and $\gmw$ and $A$ are viewed
as a time-dependent density and $1$-form on $B$, so that a dot denotes
differentiation with respect to $t$, while $D^B\gmw$ and $dA$ are the
derivatives on $B$.

Conversely for any solution of these equations on an Einstein-Weyl space $B$,
the conformal structure given by the above formula is selfdual, and $\pi$
defines a selfdual conformal submersion over $B$.
\proofof{prop} A conformal submersion certainly has the form given.  The aim
of the proof is to show that the given equations on $B$ are equivalent to the
fact that $\pi^*D^B=D^0+\omega$ on $B\times\R$.  To do this, I will work in
the (arbitrarily chosen) gauge $g=\gmw^2\conf$ and rewrite the equation
$\pi^*D^B
=D^0+\omega=D^g-\frac13(\divg^g\xi)\xi+(d^g\xi)(\xi,.)-(*d^g\xi)(\xi,.)$ using
the fact that $D^g\gmw=0$.

In the chosen gauge, $\gmw\xi=dt+A$ and so
\begin{equation*}
\gmw d^g\xi=d(\gmw\xi)=dA+dt\wedge \dot A=\gmw\xi\wedge\dot A
+dA+\dot A\wedge A.
\end{equation*}
It follows that $(d^g\xi)(\xi,.)=\dot A$ and
\begin{equation*}
(*d^g\xi)(\xi,.)={*(\xi\wedge d^g\xi)}=-{\stB(dA+\dot A\wedge A)}.
\end{equation*}
Writing $\divg^g$ in terms of $*d^g*$ readily yields
$\frac13\divg^g\xi=\dot\gmw$. Therefore:
\begin{align*}
0&=D^g\gmw=(\pi^*D^B)\gmw-\dot\gmw\gmw\xi+\dot A\gmw
+{\stB (dA+\dot A\wedge A)}\\
&=D^B\gmw+\dot \gmw dt-\dot\gmw(dt+A)+\dot A\gmw
+{\stB (dA+\dot A\wedge A)}\\
&=D^B\gmw-A\dot\gmw+\dot A\gmw+{\stB (dA+\dot A\wedge A)}.
\end{align*}
This completes the proof and also shows that $D^0=D^g-\dot\gmw\xi+\dot A$.
\end{proof}
\begin{rems} These equations may be viewed as Einstein-Weyl Bogomolny
equations for a diffeomorphism group: if the fibres are diffeomorphic to an
oriented $1$-manifold $\T$ ($S^1$ or $\R$, assuming the fibres are connected),
then $M=P\times_{\Diff(\T)}\T$ where $\Diff(\T)$ is the group of orientation
preserving diffeomorphisms of $\T$ and $P$ is the principal $\Diff(\T)$-bundle
whose fibre at $x\in B$ consists of the orientation preserving diffeomorphisms
$\T\to M_x$. Note that $P$ and $\Diff(\T)$ only enter into this formulation
infinitesimally, so one can assume that $M$ is an open subset of
$P\times_{\Diff(\T)}\T$. Choosing a gauge, i.e., a (local) section of $P$,
identifies $M$ (locally) with $B\times\T$ and the connection $1$-form and
relative length scale may be viewed as sections of $L^{-1}_B\tens\Vect(\T)$
and $T\dual B\tens\Vect(\T)$, where $\Vect(\T)$ is the Lie algebra of vector
fields on $\T$. If $t$ is a coordinate on $\T$, then writing
$\gmw=\gmw(t)\,d/dt$ and $A=A(t)\,d/dt$ shows that the equations of the above
proposition are:
\begin{equation*}
\stB (D^B\gmw+[A,\gmw])=F^A:=dA+\tfrac12[A\wedge A]
\end{equation*}
where $[.\,,.]$ denotes the Lie bracket in $\Vect(\T)$.

The classical Jones-Tod construction arises by reduction to a one dimensional
translational subgroup $S^1$ or $\R$---this point of view will be further
justified later by studying the other finite dimensional subgroups of
$\Diff(\T)$.

The gauge freedom is of course the choice of $t$ coordinate for these
monopole equations. If $\tilde t=f(t)$ for a function $f$ on $B\times\T$ with
$\dot f\neq 0$ then $\gmw(t)=\tilde\gmw(\tilde t)/\dot f$ and $A(t)=(\tilde
A(\tilde t)+df)/\dot f$. In the classical Jones-Tod correspondence there is a
preferred gauge in which to work: the constant length gauge of the conformal
vector field $K$ ($D^0=D^{|K|}$ is exact). Choosing a section $t=0$ of $M$ over
$B$ makes it into a line bundle and $M$ may be recovered from a linear
differential equation on $B$:
\begin{align*}
\conf&=\pi^*\conf\low_B+\gmw^{-2}(dt+A)^2\\ \tag*{where}\stB D^B\gmw&=dA
\quad\text{for} \quad \gmw\in\Cinf(B,L^{-1}_B).
\end{align*}
This equation for $(\gmw,A)$ is often called the (abelian) monopole equation.
\end{rems}

Now suppose that $J$ is an antiselfdual almost complex structure on $M$ which
is invariant with respect to $\xi$, i.e., $\cL^0_{\smash[b]\xi} J=0$. Let
$\Omega\low_J(X,Y)=\cip{JX,Y}$ be the conformal K{\"a}hler form of $J$ and $D$
the unique Weyl derivative such that $d^D\Omega\low_J=0$. Then it is well
known that $J$ is integrable if and only if $DJ=0$; $D$ is then called the
\emphdef{K{\"a}hler-Weyl} connection of $J$.
If $DJ=0$ then $(\cL^0_\xi J)X=J(D^0\tsum D)\low_X\xi-(D^0\tsum
D)\low_{JX}\xi$, and so $J$ is invariant iff $\xi$ is \emphdef{holomorphic} in
the sense that $(D^0\tsum D)\xi$ is complex linear.
The following generalizes a theorem of~\cite{CP2} to conformal submersions.

\addtocounter{thmcounter}{1}
\begin{thm}\label{complex} Let $M$ be an oriented conformal $4$-manifold with
a selfdual conformal submersion $\xi$ over a manifold $B$ and suppose that $J$
is an invariant antiselfdual almost complex structure on $M$.  Let
$D=D^{sd}-\kappa\xi-\tau\chi$ for basic sections $\tau$ and $\kappa$ of
$L^{-1}$, where $\chi=J\xi$. Then $DJ=0$ iff
$D^B\chi=\tau(\iden-\chi\tens\chi)+\kappa\,{\stB\chi}$ on $B$. Hence $J$ is
integrable iff $\chi$ is a shear-free geodesic congruence, in the sense that
$D^B\chi$ has the above form.
\proofof{thm} Since $J$ is invariant, $\chi$ is invariant, hence basic, since
it is horizontal. Note that $D\low_{\smash[b]\xi} J=D^{sd}_\xi
J-\tau\abrack{\xi\skwend\chi,J} =D^{sd}_\xi J=\cL^0_\xi J=0$, so it remains to
compute $D\low_XJ$ for horizontal vector fields $X$. Since
$\Omega\low_J=\xi\wedge\chi-*(\xi\wedge\chi)$ it follows that
\begin{align*}
D\low_X\Omega\low_J&=D\low_X\xi\wedge\chi+\xi\wedge D\low_X\chi
-*(D\low_X\xi\wedge\chi+\xi\wedge D\low_X\chi),\\
\tag*{where}
D\low_X\xi\wedge\chi&=D^{sd}_X\xi\wedge\chi-\kappa X\wedge\chi\\
\tag*{and} \xi\wedge D\low_X\chi&=
\xi\wedge D^{sd}_X\chi-\tau\xi\wedge\bigl(X-\cip{\chi,X}\chi\bigr).\\
\intertext{Therefore}
\xi\wedge D\low_X\chi-*(D\low_X\xi\wedge\chi)&=
\xi\wedge D^B_X\chi-\tau\xi\wedge\bigl(X-\cip{\chi,X}\chi\bigr)
+\kappa\,{*(X\wedge\chi)}.
\end{align*}
Since the right hand side is a vertical $2$-form, it follows that $D\low_XJ=0$
iff
\begin{equation*}\notag
D^B_X\chi-\cip{D^B_X\chi,\xi}=\tau(X-\cip{\chi,X}\chi)
+\kappa\,\iota\low_X\,{\stB \chi}.
\end{equation*}
To prove the final statement, suppose that $J$ is invariant and integrable
with K{\"a}hler-Weyl connection $D$.  Then $(D^0\tsum D)\xi
=-\kappa\,\iden+\frac12\tau J+\frac12(d^{D^0\tsum D}\xi)^\sd$, where
$(d^{D^0\tsum D}\xi)^\sd$ is a selfdual $2$-form and $\kappa,\tau$ are
sections of $L^{-1}$. It follows that
\begin{align*} (d^{D^0\tsum D}\xi)(\xi,.)&=
\hphantom{-}\tau\chi+(d^{D^0\tsum D}\xi)^\sd(\xi,.)\\
({*d^{D^0\tsum D}\xi})(\xi,.)&=-\tau\chi+(d^{D^0\tsum D}\xi)^\sd(\xi,.).
\end{align*}
Therefore $D^{sd}=D+\kappa\xi+\tau\chi$. It remains to check that $\kappa$ and
$\tau$ are basic. Since $DJ=0$ and $\cL^0_\xi J=0$ it follows that
$\abrack{\cL^0_\xi D\low_X,J}=0$. By the linearized Koszul formula, this
implies $\cL^0_\xi D =\frac12F^0(\xi)$ on $L^1$ and therefore
\begin{multline*}
R^D_{\xi,X}+(D^0\tsum D)\low_X
\bigl(-\kappa\,\iden+\tfrac12\tau J+\tfrac12(d^{D^0\tsum D}\xi)^\sd\bigr)\\
=\cL^0_\xi D\low_X=\tfrac12F^0(\xi,X)\iden+\tfrac12F^0(\xi)\skwend X
-\tfrac12F^0(X)\skwend\xi+\tfrac12\cip{\xi,X}F^0.
\end{multline*}
The identity and $J$ components of this formula give
$D^0\kappa+\frac12F^0(\xi)=F^D(\xi)$ and $D^0\tau=\rho^D(\xi)$ where
$\rho^D(X,Y)$ is the Ricci form of $D$, defined to be the contraction of
$R^D_{X,Y}$ with $J$. In particular $D^0_\xi\kappa=0=D^0_\xi\tau$.
\end{proof}
\begin{rem} $\tau$ and $\kappa$ are called the \emphdef{divergence} and
\emphdef{twist} of the congruence $\chi$.
\end{rem}

Assume now that $W$ is selfdual. Then so are $F^D$ and $\rho^D$ (see
e.g.~\cite{CP2}), and hence they are uniquely determined by their contractions
with $\xi$. It follows that $(M,J)$ is locally hypercomplex iff $D^0\tau=0$
iff $\tau=0$ or $D^0$ is exact and $\tau$ is constant in this gauge. This
implies that a hyperCR structure on the Einstein-Weyl quotient $B$ induces a
hypercomplex structure on $M$, generalizing a result of Gauduchon and
Tod~\cite{GT} to conformal submersions. On the other hand, $(M,J)$ is locally
scalar-flat K{\"a}hler iff $D^0\kappa+\frac12F^0(\xi)=0$, so the presence of
$F^0$ obstructs a naive generalization of LeBrun's work~\cite{LeBrun1} to this
context. This will be remedied in the next section.

I next generalize a result of Mason and Tod, which was used by
Tod~\cite{Tod3} to give a general description of selfdual Einstein metrics
with a Killing vector field.

\addtocounter{thmcounter}{1}
\begin{thm}\label{MT} Let $(M,\conf,D^{ew})$ be a selfdual Einstein-Weyl
$4$-manifold and let $\xi$ a selfdual conformal submersion with minimal Weyl
derivative $D^0$.  Then $M$ admits a canonical compatible K{\"a}hler-Weyl
structure on the open set where the antiselfdual part of
$(D^0\tsum D^{ew})\xi$ is nonzero.

More precisely, if this antiselfdual part is $\tau J$ where $J^2=-\iden$ then
$J$ is integrable, with K{\"a}hler-Weyl connection $D=D^{ew}-\tau^{-1}D^0\tau
=D^{ew}+D^0-D^\tau$, where $D^\tau$ is defined by $D^\tau\tau=0$.
\proofof{thm} It suffices to prove that $DJ=0$. Observe first that
$(D^0_X\tau)J+\tau D^{ew}_XJ$ is the antiselfdual part of
$(D^0\tsum D^{ew})\low_X(D^0\tsum D^{ew})\xi
=\cL^0_\xi D^{ew}_X-R^{ew}_{\xi,X}$,
where $R^{ew}$ is the curvature of $D^{ew}$. Expanding the curvature
and using the linearized Koszul formula to compute the Weyl-Lie derivative
gives
\begin{align*}
&(D^0_X\tau)J+\tau D^{ew}_XJ\\&\;=\bigl[\tfrac12 F^{ew}(\xi)\skwend X
+\tfrac12 F^{ew}(X)\skwend\xi+\tfrac1{12}\scal^{ew}\xi\skwend
X-D^0\kappa\skwend X+\tfrac12\xi\skwend F^0(X)\bigr]^\asd\\
&\;=\bigl[(F^{ew}(\xi)+\tfrac1{12}\scal^{ew}\xi-D^0\kappa
-\tfrac12F^0(\xi))\skwend X\bigr]^\asd,
\end{align*}
where $[...]^\asd$ denotes the antiselfdual part, $\kappa$ is minus the
identity component of $(D^0\tsum D^{ew})\xi$, $\scal^{ew}$ is the scalar
curvature of $D^{ew}$, $F^{ew}$ is the Faraday curvature of $D^{ew}$, and I
have used the fact that $F^0$ and $F^{ew}$ are selfdual. The precise form of
this expression is now not important: it suffices to observe that it is of the
form $[J\alpha\skwend X]^\asd$ for some $1$-form $\alpha$. Since
${*(J\alpha\skwend X)}=-(\alpha\skwend JX+\alpha(X)J)$, it follows that
$J\alpha\skwend X-{*(J\alpha\skwend X)}=\abrack{\alpha\skwend X,J}
+\alpha(X)J$ and the commutator term is orthogonal to $J$. Since $D^{ew}J$ is
also orthogonal to $J$, $D^0_X\tau=\alpha(X)$ and $\tau
D^{ew}_XJ=\abrack{\alpha\skwend X,J}=\abrack{D^0\tau\skwend X,J}$, i.e., $J$ is parallel
with respect to $D^{ew}-\tau^{-1}D^0\tau$.
\end{proof}

The theorems of this section reduce to known results when $D^0$ is exact
(i.e., $\xi=K/|K|$ for some conformal vector field $K$), but they have one
disadvantage over the results they generalize.  Namely, the abelian monopole
equation on $B$ arising in the classical Jones-Tod correspondence becomes a
nonlinear evolution equation, which is much harder to solve.
This difficulty has already been encountered in another special case of the
above theorems: the case that $M$ is the selfdual Einstein metric locally
``filling in'' $B$ via Hitchin's version of LeBrun's $\cH$-space
construction~\cite{Hitchin3,LeBrun0}. This beautiful construction of a
selfdual Einstein metric from a real analytic conformal $3$-manifold (which is
taken to be Einstein-Weyl in Hitchin's construction) is defined twistorially,
making it difficult to carry out in practice.

Furthermore conformal submersions themselves are hard to find, because the
equation for conformal submersions, unlike the conformal Killing equation, is
nonlinear.  Hence, for the theorems of this section to be interesting, it is
essential to find new situations in which the inverse construction can be
carried out and examples can be found. This will be done in the next two
sections. In the final section, a direct version of the Hitchin-LeBrun
construction will be obtained.

\section{Affine conformal submersions}\label{acs}

An affine structure on a submersion $\pi\colon M\to B$ is a flat torsion-free
connection on each fibre. This identifies $M$, at least locally, with an
affine bundle modelled on the vector bundle whose fibre at each point of $B$
is the space of parallel vector fields on the corresponding fibre of $M$. For
conformal submersions with oriented one dimensional fibres, the vertical
bundle $\Vb$ of $M$ is isomorphic to $L^1$ and so the vertical part
$D\low_\xi$ of any Weyl derivative $D$ defines an affine structure on $M$.

There are many choices of affine structure on $M$, but such a choice is only
helpful if the conformal structure on $M$ is affine; that is, in terms of
Proposition~\ref{confsub}, the nonlinear connection on $M$ is an affine
connection, and the relative length scale is affine along the fibres. If such
a ``good'' affine structure can be found, I will say that $\pi$ is an
\emphdef{affine conformal submersion}.

More precisely, a nonlinear connection $\Hb$ induces a linearized connection
on the infinite dimensional vector space of vertical vector fields defined by
$\cD_X U=[\XH,U]$, where $U$ is a vertical vector field, $X$ is a vector field
on $B$ and $\XH$ is its horizontal lift, so that $[\XH,U]$ is vertical. $\Hb$
is affine iff the parallel vertical vector fields on each fibre are preserved
by $\cD$; this then induces the linearized connection on the model vector
bundle. Affine connections form an affine space modelled on $1$-forms on $B$
with values in the affine vector fields on $\Vb$. Similarly, the relative
length scale $\gmw\colon\pi^*L^1_B\to\Vb$ is affine iff it maps basic
densities to affine vector fields, in which case it may be viewed as a
$(-1)$-density on $B$ with values in the affine vector fields.

The above approach and the next proposition arose from discussions with Paul
Gauduchon in a joint effort to understand affine conformal submersions.

\begin{prop} Let $M$ be a conformal manifold and let $\xi$ be a
conformal submersion over $B$ with minimal Weyl derivative $D^0$. Define an
affine structure $D\low_{\smash[b]\xi}$ on $\pi\colon M\to B$ by the Weyl
derivative $D=D^0+\lambda\xi$ where $\lambda$ is a section of $L^{-1}$.
\begin{enumerate}
\item The connection $\Hb$ on $M\to B$ is
affine with respect to $D\low_{\smash[b]\xi}$ iff $F^D(\xi)=0$.
\item The relative length scale is affine with respect to $D\low_\xi$
iff $\lambda$ is basic.
\end{enumerate}
Hence the conformal submersion is affine if $D^0\lambda=F^0(\xi)$.
\proofof{prop} The $D\low_{\smash[b]\xi}$-parallel vertical vector fields are
defined by identifying $\Vb$ with $L^1$ using $\xi$.  Hence the linearized
connection may be defined on $\mu\in\Cinf(M,L^1)$ by
\begin{align*}
(\cD_X\mu)\xi&=[\XH,\mu\xi]=D^0_\XH(\mu\xi)-\mu D^0_\xi\XH=(D^0_\XH\mu)\xi-
\mu(D^0_\XH\xi-D^0_\xi\XH)\\
&=(D^0_\XH\mu)\xi+\cL^0_\xi\XH=(D\low_\XH\mu)\xi
\end{align*}
since $\XH$ is invariant and $D-D^0$ is vertical. Hence $\Hb$ is affine iff
$D\low_{\smash[b]\xi}(D_\XH\mu)=0$ for all $\mu$ with
$D\low_{\smash[b]\xi}\mu=0$.  Since $[U,\XH]$ is vertical for $U$ vertical,
$D\low_{\smash[b]\xi}(D_\XH\mu)=F^D(\xi,\XH)\mu$ and so $\Hb$ is affine iff
$F^D(\xi)=0$.

The relative length scale is section $\gmw$ of $\pi^*L^{-1}_B\tens\Vb$. This
is affine iff its vertical derivative with respect to the affine structure, as
a section of $\pi^*L^{-1}_B\tens\Vb\dual\tens\Vb\cong\pi^*L^{-1}_B$, is basic.
Identifying $\pi^*L^1_B$ and $\Vb$ with $L^1$ identifies $\gmw$ with the
identity map in $L^{-1}\tens L^1$ but its vertical derivative must be computed
with respect to the covariant derivative $D^0\tsum D$ and so $\gmw$ is affine
iff $0=D^0_\xi(D^0\tens D)\low_\xi\iden=D^0_\xi\lambda$.

Now observe that $F^D(\xi)=F^0(\xi)+d(\lambda\xi)(\xi)=F^0(\xi)
+(D^0_\xi\lambda)\xi-D^0\lambda$.
\end{proof}
It follows from this that there is an obstruction to the existence of a good
affine structure for a conformal submersion: since
$-F^0\lambda=d^0(F^0(\xi))=\cL^0_\xi F^0$, the Weyl-Lie derivative of $F^0$
must be a multiple of $F^0$. If a good affine structure exists, it is
essentially unique: any two differ by a section $\mu$ of $L^{-1}$ with
$D^0\mu=0$ which implies that the affine structures are equal or $D^0$ is
exact and $\mu$ is constant.

I now return to four dimensions and the Jones-Tod construction.

\addtocounter{thmcounter}{1}
\begin{thm}\label{affine} Let $(M,\conf)$ be a selfdual conformal $4$-manifold
with a selfdual affine conformal submersion $\pi\colon M\to B$ over an
Einstein-Weyl space $(B,\conf\low_B,D^B)$.  Then with respect to an arbitrary
affine coordinate $t$ on $M\to B$, the conformal structure on $M$ is
\begin{align*}
\conf=\pi^*\conf\low_B &+ (t\gmw_1 +\gmw_0)^{-2}(dt+tA_1+A_0)^2\\
\tag*{\textit{where}} \stB  D^B\gmw_1&=dA_1\\
\tag*{\textit{and}}\stB (D^B\gmw_0+A_1\gmw_0-A_0\gmw_1)&=dA_0+A_1\wedge A_0
\end{align*}
for some $\gmw_0,\gmw_1\in \Cinf(B,L^{-1}_B)$ and $A_0,A_1\in\Cinf(B,T\dual
B)$. Conversely for any solution of these affine monopole equations on an
Einstein-Weyl space $B$, the conformal structure given by the above formula is
selfdual, and the above decomposition defines a selfdual affine conformal
submersion over $B$.
\proofof{thm} An affine conformal submersion certainly has the form given.
It remains to apply this Ansatz to the equations of Proposition~\ref{evolve},
by writing $\gmw=t\gmw_1+\gmw_0$ and $A=tA_1+A_0$. Now
\begin{align*}
D^B\gmw+\dot A\gmw-A\dot\gmw&=tD^B\gmw_1+D^B\gmw_0+\gmw_0A_1-\gmw_1A_0\\
dA+\dot A\wedge A&=t dA_1+dA_0+A_1\wedge A_0
\end{align*}
and the linear and constant terms (in $t$) of these equations
prove the result.
\end{proof}
\begin{rems} Note that $D^0=D^g-\gmw_1\xi+A_1$ and so
$F^0(\xi)=-d(\gmw_1\xi)(\xi)=D^0\gmw_1-(D^0_\xi\gmw_1)\xi=D^0\gmw_1$ since
$\gmw_1$ is basic. Hence $\lambda=\gmw_1$ is the solution of
$D^0\lambda=F^0(\xi)$: the affine structure is $D^g_\xi=D^0_\xi+\gmw_1$.

The freedom in the choice of affine coordinate $t$ gives a gauge freedom for
the affine monopole equations. If $\tilde t=at+b$ for basic functions $a,b$,
write $\gmw_1=\tilde\gmw_1$, $\gmw_0= a^{-1}(\tilde\gmw_0+b\tilde\gmw_1)$,
$A_1=\tilde A_1+a^{-1}da$ and $A_0=a^{-1}(\tilde A_0+b\tilde A_1+db)$. One
immediately verifies, by substituting into the affine monopole equations, that
$(\tilde\gmw,\tilde A)$ is a solution if $(\gmw,A)$ is. Note that
$t\mu_g=t\gmw^{-1}$ is a well defined section of $L^1$ up to translation by a
basic section of $L^1$: it may be fixed by choosing a section of $M\to B$. The
induced exact Weyl derivative is
$D^{t\mu_g}=D^g-t^{-1}dt=D^0-t^{-1}\gmw_0\xi+t^{-1}A_0$.

The equations show that the affine Jones-Tod correspondence reduces to the
classical case in two ways. Firstly the linear part of the affine monopole
equation is an abelian monopole equation: if $(\gmw_0,A_0)$ is zero, the
affine bundle $M$ is isomorphic to the model vector bundle.  Secondly if the
linear part $(\gmw_1,A_1)$ vanishes, the translational part of the of the
affine monopole equation is an abelian monopole equation: the model vector
bundle is trivial, and so $M$ is a principal $\R$-bundle. On the other hand if
the solution $(\gmw_1,A_1)$ is nontrivial, then it gives a linearization of
$M$, i.e., a selfdual space with an affine conformal submersion is affinely
modelled on a selfdual space with a conformal vector field.
\end{rems}

This theorem gives a new method for constructing selfdual spaces from linear
equations, since the second affine monopole equation is linear once a solution
of the first equation is chosen.  In particular, this method gives all
scalar-flat K{\"a}hler metrics with a holomorphic selfdual conformal submersion
including all hyperK{\"a}hler metrics admitting such a submersion.

\addtocounter{thmcounter}{1}
\begin{thm}\label{SFK} Let $(M,g)$ be a four dimensional scalar-flat K{\"a}hler
metric \textup(with antiselfdual complex structure, so that $(M,\conf)$ is
selfdual\textup) admitting a holomorphic selfdual conformal submersion. Then
the conformal submersion is affine over an Einstein-Weyl space with a
shear-free geodesic congruence $\chi$, where the linear part of the affine
monopole is given by $\gmw_1=-2\kappa$ and $\kappa$ is the twist of the
congruence $\chi$. All such metrics are locally of the form
\begin{equation*}
g=(\rho-2\mu_t^{-1}\kappa)\conf\low_B
+\frac{\bigl(D^B(\mu_t^{-1})+2\tau\chi\mu_t^{-1}+\Phi\bigr)^2}
{\rho-2\mu_t^{-1}\kappa},
\end{equation*}
where $\rho\in\Cinf(B,L^{-2}_B)$, $\Phi\in\Cinf(B,L^{-1}_BT\dual B)$,
$\mu_t^{-1}$ is a section of $L^{-1}$ increasing along the fibres, and $\tau$
is the twist of the congruence $\chi$. Conversely, for any Einstein-Weyl space
$(\conf\low_B,D^B)$ and shear-free geodesic congruence $\chi$ with twist
$\kappa$ and divergence $\tau$, this metric is scalar-flat K{\"a}hler iff
$(\rho,\Phi)$ satisfy the linear differential equation
\begin{equation*}
\stB(D^B\rho+2\tau\chi\rho+2\kappa\Phi)=d^B\Phi+2\tau\chi\wedge\Phi.
\end{equation*}
Furthermore the metric is hyperK{\"a}hler iff $D^0\tau=0$, i.e., iff
$\tau=0$ or $D^0$ is exact and $\tau$ is constant in this gauge.
\proofof{thm} By Theorem II, $D^g=D^{sd}-\kappa\xi-\tau\chi$ where
$D^B\chi=\tau(\iden-\chi\tens\chi)+\kappa\stB \chi$, and furthermore,
$D^0\kappa+\frac12F^0(\xi)=0$, since $F^{D^g}=0$. Hence $D^0_\xi-2\kappa$
defines an affine structure on $M$ making the submersion affine and
$\gmw_1=-2\kappa$.

The analysis of shear-free geodesic congruences in~\cite{CP2}
shows that $\stB D^B\kappa=\frac12 F^B-d(\tau\chi)$, and hence,
choosing any gauge $\muB$, one can take $A_1=-\omB+2\tau\chi$.
The translational part of the monopole equation is therefore:
\begin{equation*}
\stB (D^B\gmw_0-\omB\gmw_0+2\tau\chi\gmw_0+2\kappa A_0)
=dA_0+(-\omB+2\tau\chi)\wedge A_0.
\end{equation*}
Now the Levi-Civita connection of the affine gauge is
$D^0+\gmw_1\xi-A_1=D^0-2\kappa\xi-2\tau\chi+\omB$, whereas the Levi-Civita
connection of the $\muB$-gauge is $D^B-\omB$. The barycentre of these
is $D^{sd}-\kappa\xi-\tau\chi$, which is the Levi-Civita connection of the
scalar-flat K{\"a}hler metric. This identifies $g$ within the conformal class,
and putting $\rho=\mu_B^{-1}\gmw_0$, $\Phi=\mu_B^{-1}A_0$ and
$\mu_t^{-1}=t\mu_B^{-1}$ completes the proof.
\end{proof}
When $\kappa=0$ this is LeBrun's construction of scalar-flat K{\"a}hler metrics
with Killing fields~\cite{LeBrun1}. On the other hand, when $(\gmw_0,A_0)=0$,
this theorem reduces to the construction of scalar-flat K{\"a}hler metrics with
homothetic vector fields~\cite{CP2}, including, as a special case, the
hyperK{\"a}hler metrics of~\cite{GT}.

I end this section with some nontrivial examples. In~\cite{CT}, the following
Einstein-Weyl structures were found from solutions of the $\SU(\infty)$ Toda
field equation.
\begin{equation*}
g\low_B=(z+h)(z+\overline h)g_{S^2}+dz^2,\qquad
\omB =-\frac{2z+h+\overline h}{(z+h)(z+\overline h)}dz,
\end{equation*}
where $h$ is a holomorphic function on an open subset of $S^2$ and
$D^B=D^g+\omega$. Note that the weightless unit vector field dual to $dz$
generates a shear-free geodesic congruence with vanishing twist
($\tau\neq0,\kappa=0$). These spaces also admit shear-free geodesic
congruences with vanishing divergence ($\tau=0,\kappa\neq0$), i.e., they are
hyperCR, and they are called the \emphdef{hyperCR-Toda spaces}.

Applying the classical Jones-Tod construction to these spaces gives conformal
structures of the form $\conf=\pi^*\conf\low_B+\gmw^{-2}(\beta+\gmv\,dz)^2$,
where $\beta=dt+\theta$ for a $1$-form $\theta$ on $B$ orthogonal to $dz$, and
$\stB D^B\gmw=dA$ with $A=\theta+\gmv\,dz$. 

These conformal structures admit a compatible scalar-flat K{\"a}hler metric and
also a compatible hypercomplex structure. The conformal vector field
$\partial/\partial t$ is a Killing field of the scalar-flat K{\"a}hler metric
and triholomorphic with respect to the hypercomplex structure.

For certain solutions of the abelian monopole equation, $\partial/\partial z$
defines a conformal submersion. To see this, write
$\conf=\eps_0^2+\eps_1^2+\eps_2^2+\eps_3^2$ where $\eps_0$ and $\eps_3$ are
the weightless unit 1-forms corresponding to $\gmw\,dz$ and
$\beta+\gmv\,dz$. The weightless unit $1$-form dual to $\partial/\partial z$
is $\xi=(\gmw\eps_0+\gmv\eps_3)/\sqrt{\gmw^2+\gmv^2}$ and so
$\eps_0^2+\eps_3^2-\xi^2=(\gmv\eps_0-\gmw\eps_3)^2/(\gmw^2+\gmv^2)=
\gmw^2\beta^2/(\gmw^2+\gmv^2)$. Hence if $\beta$ and $(\gmw^2+\gmv^2)|z+h|^2$
are independent of $z$, then $\partial/\partial z$ will define a conformal
submersion with quotient $(\gmw^2+\gmv^2)|z+h|^2g_{S^2}+\beta^2$.

Now Ian Strachan has pointed out~\cite{CT} that for any holomorphic function
$f$,
\begin{equation*}\notag
\gmw =\frac12\left(\frac{f}{z+h}+\frac{\overline f}{z+\overline h}\right),
\quad
\gmv =\frac1{2i}\left(\frac{f}{z+h}-\frac{\overline f}{z+\overline h}\right),
\quad d\beta=\tfrac12(f+\overline f)\vol_{S^2}
\end{equation*}
defines a solution of the monopole equation. Clearly $\beta$ and
$(\gmw^2+\gmv^2)|z+h|^2$ are independent of $z$ and so $\partial/\partial z$
defines a conformal submersion. Explicitly, $\conf$ has a compatible
metric
\begin{equation*}
g=\biggl(\frac{z+h}{2f}+\frac{z+\overline h}{2\overline f}\biggr)^2
\big(|f|^2g_{S^2}+\beta^2\bigr)+\left[dz+i
\biggl(\frac{z+h}{2f}-\frac{z+\overline h}{2\overline f}\biggr)\beta\right]^2
\end{equation*}
and so this submersion is obviously affine, with affine coordinate $z$. The
quotient conformal $3$-manifold $\tilde B$ admits an Einstein-Weyl structure:
\begin{equation*}
g\low_{\tilde B}=|f|^2g_{S^2}+\beta^2,\qquad
\omega\low_{\!\tilde B}=\frac i2\biggl(\frac1f-\frac1{\overline f}\biggr)\beta.
\end{equation*}
These are the Einstein-Weyl spaces with \emphdef{geodesic symmetry} described
in~\cite{CP2}. One easily checks that $g$ is given by a solution of the affine
monopole equations with $\gmw_1=-2\kappa_s$ where $\kappa_s$ is the twist of
the geodesic symmetry on $\tilde B$, i.e., $d\beta=2\kappa_s{\stB\beta}$.
Hence these scalar-flat K{\"a}hler metrics with compatible hypercomplex structures
could have been constructed directly as selfdual affine conformal submersions
over the Einstein-Weyl spaces with geodesic symmetry. When $f=ah+b$ for
$a,b\in\R$, these metrics are conformally Einstein~\cite{CT} and will feature
again in the final section.

\section{Projective conformal submersions}\label{pcs}

A natural generalization of an affine conformal submersion is a projective
conformal submersion. A projective structure on a $1$-manifold is a second
order linear differential operator from $L^{1/2}$ to $L^{-3/2}$ which has no
first order term with respect to any Weyl derivative, and the same definition
may be applied fibrewise to a congruence $\xi$. Hence any Weyl derivative $D$
induces a projective structure $\mu\mapsto
D\low_{\smash[b]\xi}(D\low_{\smash[b]\xi}\mu)$. Note that
$(D+\gamma)\low_\xi\bigl((D+\gamma)\low_{\smash[b]\xi}\mu\bigr)=
(D\low_{\smash[b]\xi}-\frac12\gamma(\xi))
(D\low_{\smash[b]\xi}\mu+\frac12\gamma(\xi)\mu)
=D\low_{\smash[b]\xi}(D\low_{\smash[b]\xi}\mu)
+\frac12D\low_{\smash[b]\xi}(\gamma(\xi))\,\mu -\frac14(\gamma(\xi))^2\,\mu$,
verifying that the condition of vanishing first order term is independent of
the Weyl derivative.

A conformal submersion $\pi$ will be called \emphdef{projective} iff there is
a projective structure on $\pi$ such that the connection $\Hb$ is projective
and the relative length scale $\gmw$ takes values in the projective vector
fields: recall that these are characterized as being quadratic in any
projective coordinate.
\begin{rem}
A curve in a conformal manifold, with weightless unit tangent $\xi$ has a
canonical projective structure given by $(D\low_{\smash[b]\xi})^2
+\frac12r^D(\xi,\xi)+\frac14|D\low_{\smash[b]\xi}\xi|^2$,
and is called a conformal geodesic if
$D\low_{\smash[b]\xi}(D\low_{\smash[b]\xi}\xi)+|D\low_{\smash[b]\xi}\xi|^2\xi
-r^D(\xi)+r^D(\xi,\xi)\xi=0$; these expressions are independent of the Weyl
derivative $D$. However, if $\xi$ is a projective conformal submersion,
there is no reason for the projective structure to equal the canonical one,
nor will the fibres be conformal geodesics in general.
\end{rem}
\begin{prop} Let $M$ be a conformal manifold and let $\xi$ be a conformal
submersion over $B$ with minimal Weyl derivative $D^0$.  Then the conformal
submersion is projective with respect to the projective structure
$(D^0_\xi)^2+\tfrac12\rho$, for a section $\rho$ of $L^{-2}$, iff
$D^0\rho=\cL^0_\xi F^0(\xi)$.

\proofof{prop} Write the projective structure as $(D_\xi)^2$ where
$D_\xi=D^0_\xi+\lambda$ is a compatible affine structure, so that
$\rho=D^0_\xi\lambda-\frac12\lambda^2$.  Then the connection
$\Hb$ is projective iff it maps $D_\xi$-parallel vertical vector fields to
$D_\xi$-affine vertical vector fields, i.e., iff $D_\xi\mu=0\implies
\partial_\xi(D_\xi(D_X\mu))=0$ for basic vector fields $X$. This condition
reduces easily to $(\cL^0_\xi-\lambda)(F^0(\xi,X)-D^0_X\lambda)=0$.

The relative length scale is projective iff
$0=(D^0\tsum D)\low_\xi(D^0_\xi\lambda)=(\cL^0_\xi-\lambda)(D^0_\xi\lambda)$.
Hence the conformal submersion is projective iff
\begin{align*}
0&=(\cL^0_\xi-\lambda)(F^0(\xi)-D^0\lambda)\\
&=(\cL^0_\xi-\lambda)F^0(\xi)-(\cL^0_\xi D^0)\lambda-D^0(\cL^0_\xi\lambda)
+\lambda D^0\lambda\\
&=\cL^0_\xi F^0(\xi)-D^0(D^0_\xi\lambda-\tfrac12\lambda^2)
\end{align*}
since $\cL^0_\xi D^0=-F^0(\xi)$ on $L^{-1}$.
\end{proof}
There is still an obstruction to solving this, since it implies
$-2F^0(\xi)\rho=\cL^0_\xi\cL^0_\xi F^0(\xi)$.

\addtocounter{thmcounter}{1}
\begin{thm}\label{projective} Let $(M,\conf)$ be a selfdual conformal
$4$-manifold with a selfdual projective conformal submersion $\pi\colon M\to
B$ over an Einstein-Weyl space $(B,\conf\low_B,D^B)$. Then with respect to an
arbitrary projective coordinate $t$ on $M\to B$, the conformal structure on
$M$ is
\begin{align*}
\conf=\pi^*\conf\low_B
+(t^2\gmw_2+t\gmw_1+\gmw_0)^{-2}&(dt+t^2A_2+tA_1+A_0)^2,\\
\tag*{\textit{where}}
\stB \bigl(D^B\gmw_2+A_2\gmw_1-A_1\gmw_2\bigr)
&=\hphantom{\tfrac12}dA_2+A_2\wedge A_1,\\
\stB \bigl(\tfrac12D^B\gmw_1+A_2\gmw_0-A_0\gmw_2\bigr)
&=\tfrac12dA_1+A_2\wedge A_0,\\
\tag*{\textit{and}}
\stB \bigl(D^B\gmw_0+A_1\gmw_0-A_0\gmw_1\bigr)
&=\hphantom{\tfrac12}dA_0+A_1\wedge A_0,
\end{align*}
for some $\gmw_0,\gmw_1,\gmw_2\in \Cinf(B,L^{-1}_B)$ and
$A_0,A_1,A_2\in\Cinf(B,T\dual B)$. Conversely for any solution of these
projective monopole equations on an Einstein-Weyl space $B$, the conformal
structure given by the above formula is selfdual, and the above decomposition
defines a selfdual projective conformal submersion over $B$.
\proofof{thm} As in the proof of Theorem IV, this amounts to computing the
equations of Proposition~\ref{evolve}, now with
$\gmw=t^2\gmw_2+t\gmw_1+\gmw_0$ and $A=t^2A_2+tA_1+A_0$. This leads to
the quadratic expressions
\begin{align*}
D^B\gmw+\dot A\gmw-A\dot\gmw&=t^2\bigl(D^B\gmw_2+A_2\gmw_1-A_1\gmw_2\bigr)\\
&\qquad+2t\bigl(\tfrac12D^B\gmw_1+A_2\gmw_0-A_0\gmw_2\bigr)
+D^B\gmw_0+A_1\gmw_0-A_0\gmw_1\\
dA+\dot A\wedge A&=t^2(dA_2+A_2\wedge A_1)+2t(\tfrac12 dA_1+A_2\wedge A_0)
+dA_0+A_1\wedge A_0
\end{align*}
and equating coefficients (in $t$) of the resulting equations completes the
proof.
\end{proof}
Note that $D^g_\xi=D^0_\xi+2t\gmw_2+\gmw_1$ and
$D^{t\mu_g}=D^0+(t\gmw_2-t^{-1}\gmw_0)\xi-tA_2+t^{-1}A_0$.
The projective structure $(D^0_\xi)^2+\rho$ is given by
$\rho=\gmw_0\gmw_2-\frac14\gmw_1^2$.

The equations arising in this theorem may be identified as $\SL(2,\R)$
Einstein-Weyl Bogomolny equations: writing
\begin{equation*}
\gmw=\begin{pmatrix}\tfrac12\gmw_1& \gmw_0\\ -\gmw_2&
-\tfrac12\gmw_1\end{pmatrix},
\qquad A=\begin{pmatrix}\tfrac12A_1& A_0\\ -A_2&
-\tfrac12A_1\end{pmatrix},
\end{equation*}
yields an $\lie{sl}(2,\R)$-valued density and connection $1$-form. The
equations of the above Theorem now become $\stB(D^B+\ad A)\gmw=F^A:=dA+A\wedge
A$.

\section{Twistor theory of conformal submersions}\label{twistor}

The constructions discussed so far have a natural interpretation on the
\emphdef{twistor space} $Z$ of $M$. This is a complex manifold fibering over
$M$ whose fibres are the antiselfdual complex structures on each tangent space
of $M$ (see~\cite{AHS,Besse}). The antipodal map on each fibre defines a real
structure (antiholomorphic involution) $\sigma$ on $Z$, so the fibres of $Z$
are real, i.e., $\sigma$-invariant. Each fibre $Z_x$ has normal bundle
$N^{(x)}\cong\cO(1)\dsum\cO(1)$ and so the fibres are precisely the real lines
amongst the ``twistor lines'', which are the holomorphic deformations of a
typical fibre. Hence real holomorphic sections of $N^{(x)}$ over $Z_x$ are
constant as maps from $Z_x$ to $T_xM$.

Each Weyl derivative $D$ induces a connection on $\pi_Z\colon Z\to M$
and hence a projection $v^D\colon TZ\to VZ$ onto the vertical bundle of $Z$.
Under a change of Weyl derivative,
$v^{D+\gamma}(U)=v^D(U)+\abrack{\gamma\skwend d\pi_Z(U),J}$ for $U\in T_JZ$.
If $K$ is any vector field on $M$ and $J\in Z_x$, then the commutator
$\abrack{DK-\frac12\cL_K\conf,J}$ is a skew endomorphism of $T_x M$
anticommuting with $J$ (since $\cL_K\conf=2\sym_0 DK$) and hence an element of
$V_JZ=T_J(Z_x)$. The lift $K^\C$ of $K$ to $Z$ defined by
$v^D(K^\C)=\abrack{DK-\frac12\cL_K\conf,J}$ is easily seen to be independent
of the choice of $D$, and is a holomorphic vector field iff $K$ is a conformal
vector field, in which case $v^D(K^\C)=\abrack{DK,J}$ (see~\cite{Gauduchon4}).
This generalizes to congruences.

\begin{prop}
Congruences $\xi$ on $M$ \textup(up to a sign\textup) are in bijective
correspondence with complex line subbundles $L^\xi$ of $TZ$ which are
$\sigma$-invariant and transverse to the real twistor lines.

$L^\xi$ is a holomorphic subbundle of $TZ$ iff $\xi$ is a selfdual conformal
submersion. The holomorphic structure on this line bundle corresponds, under
the Ward correspondence, to the minimal Weyl derivative $D^0$.

\proofof{prop}[Sketch proof] Given $\xi$, let $L^\xi$ be the complex span of
those vectors $U$ in $T_JZ$ with $d\pi^Z(U)=\mu\xi$ and $v^D(U)=\mu
\abrack{(D^0\tsum D)\xi-\frac12\cL^0_\xi\conf,J}$ for some element $\mu$ of
$L^1_{\pi(J)}$ (recall that $\frac12\cL^0_\xi\conf=\sym_0(D^0\tsum D)\xi$).
This line subbundle $L^\xi$ of $TZ$ naturally isomorphic to $\pi^*L^1\tens\C$
and is clearly transverse to the real twistor lines.  Conversely, a
$\sigma$-invariant complex line subbundle $L^\xi$ transverse to a real twistor
line $Z_x$ must be degree $0$ and the real sections define a congruence $\xi$
up to sign, and identify $L^\xi$ with $\pi^*L^1\tens\C$.

By the Ward correspondence, a $\overline\partial$-operator on $L^\xi$
corresponds to a Weyl derivative $D^0$ on $M$, and $L^\xi$ holomorphic if
$D^0$ is selfdual~\cite{Gauduchon4}. If so, the inclusion of $L^\xi$ into $TZ$
may be viewed as a section $\xi^\C$ of the holomorphic bundle
$(L^\xi)^{-1}\tens TZ$ and one finds that $\xi^\C$ is holomorphic iff
$\cL^0_\xi\conf=0$.
\end{proof}

\noindent This gives a twistorial explanation for the theorems of
section~\ref{JTC}.

\smallbreak\noindent 1.\ The distribution $L^\xi$ on $Z$ given by a selfdual
conformal submersion $\xi$ integrates a holomorphic foliation with one
dimensional leaves.  Since $L^\xi$ is trivial on each twistor line, the
twistor lines map to rational curves (called ``minitwistor lines'') with
normal bundle $\cO(2)$ in the (local) quotient space $\cS$, which is the
``minitwistor space'' that gives rise to the Einstein-Weyl structure on
$M/\xi$~\cite{Hitchin3}.

\smallbreak\noindent 2.\ The condition that $\xi$ is holomorphic with respect
to an antiselfdual complex structure $J$ is simply the condition that the
image of $J$ (as a section of $Z$) is a union of leaves of the foliation
determined by $\xi$.  This divisor $\cD$ in $Z$ therefore descends to a
divisor $\cC$ in the minitwistor space, which in turn determines a shear-free
geodesic congruence~\cite{CP2}. The correspondence between hypercomplex and
hyperCR spaces follows from the fact that $[\cD-\overline\cD]$ is trivial if
$[\cC-\overline\cC]$ is trivial: if so this gives a map from $\cS$ to $\CP1$
and hence from $Z$ to $\CP1$. On the other hand, the correspondence between
scalar-flat K{\"a}hler metrics and Toda Einstein-Weyl spaces has a more subtle
generalization because the canonical bundle $K_Z$ is no longer the pullback of
$K_S$ and so $[\cD-\overline\cD]K_Z^{1/2}$ is not the pullback of
$[\cC-\overline\cC]K_\cS^{1/2}$.

\smallbreak\noindent 3.\ Compatible Einstein-Weyl structures $D$ on selfdual
conformal manifolds are Einstein or locally hypercomplex~\cite{DMJC2,PS1} and
correspond to holomorphic rank two distributions $H^D=\ker v^D$ on $Z$. The
twisted $1$-form $\theta\in H^0\bigl(Z,(L^D)^{-1}K_Z^{-1/2}T\dual Z\bigr)$
defining $H^D$ can be contracted with $\xi^\C\in
H^0\bigl(Z,(L^\xi)^{-1}TZ\bigr)$ to give $\theta(\xi^\C)$, a holomorphic
section of $(L^D)^{-1}(L^\xi)^{-1}K_Z^{-1/2}$, which has degree two on each
twistor line. If this section is not identically zero, then the corresponding
divisor gives rise to a complex structure. On the other hand, if the section
is identically zero, then $L^\xi$ is a subbundle of $H^D$ and so
$\abrack{(D^0\tsum D)\xi,J}=0$ for all antiselfdual almost complex structures
$J$.  \smallbreak

The affine and projective cases may also be interpreted twistorially: the
holomorphic bundle $Z\to\cS$ is an affine or projective line bundle. Such
bundles arise as affine subspaces or projectivizations of rank two vector
bundles trivial on minitwistor lines. Since these correspond to Einstein-Weyl
monopoles via a generalized Hitchin-Ward correspondence~\cite{Hitchin3a}, it
is no surprise that the equations are affine and projective monopole
equations.

When the Einstein-Weyl structure is Einstein with nonzero scalar curvature,
$H^D$ is a contact distribution, whereas when it is locally hypercomplex,
$H^D$ is integrable. Therefore, if the skew symmetric part of $(D^0\tsum
D)\xi$ is selfdual (i.e., if $L^\xi$ is a subbundle of $H^D$), I will say that
$\xi$ is \emphdef{Legendrian} or \emphdef{triholomorphic} in the case that the
scalar curvature is nonzero or zero respectively.

In the Legendrian case, the leaves of the foliation of $Z$ given by $\xi$ are
Legendrian curves in a contact manifold. Let $\cS$ be the (local) quotient and
let $Z^{(y)}$ be the leaf corresponding to a point $y\in\cS$. Then at each
point $J$ of $Z^{(y)}$, the contact distribution projects onto a one
dimensional subbundle of $T_y\cS$, giving a holomorphic map of complex curves
from $Z^{(y)}$ to $P(T_y\cS)$. This map cannot be constant, as the contact
distribution is non-integrable, so it is locally an isomorphism. By its very
definition, this isomorphism identifies the contact distribution on $Z$ with
the canonical contact distribution on $P(T\cS)\cong P(T\dual\cS)$, where a
line in $T_y\cS$ is identified with its annihilator in
$T\dual_y\cS$. Therefore, $Z$ can be locally
identified (in a neighbourhood of any twistor line) with the projectivized
cotangent bundle of $\cS$. This is Hitchin's construction of the selfdual
Einstein metric (with nonzero scalar curvature) ``filling in'' a
$3$-dimensional Einstein-Weyl space~\cite{Hitchin3}. LeBrun~\cite{LeBrun0} has
given such a construction for any real analytic conformal $3$-manifold $B$,
and Hitchin observes that the choice of a compatible Einstein-Weyl structure
on $B$ (if one exists) equips the selfdual Einstein metric with a conformal
submersion onto $B$. The discussion here \emph{characterizes} the
conformal submersions arising in this way.

\addtocounter{thmcounter}{1}
\begin{thm}\label{hitchin} Let $M$ be a selfdual Einstein manifold with nonzero
scalar curvature. Then $M$ arises from Hitchin's construction iff it admits a
Legendrian selfdual conformal submersion.
\end{thm}
The scalar curvature is usually taken to be negative: it is in such a real
slice that the original conformal $3$-manifold appears as a conformal infinity.

When the Einstein-Weyl space is hyperCR, its minitwistor space fibres over
$\CP1$ and the vertical bundle of this fibration is transverse to the
minitwistor lines. This line subbundle of $T\cS$ defines a section $X$ of
$P(T\cS)$ which does not intersect the lifted minitwistor lines, and hence
does not intersect nearby twistor lines. On removing this section
$P(T\cS)\setdif X(\cS)$ is an affine bundle over $\cS$ (and is still a twistor
space). Hence one can expect to carry out the Hitchin-LeBrun construction
explicitly in this case, using the affine monopole equations.

In general, $P(T\cS)=P(K_{\smash[b]\cS}^{1/2}T\cS)$ is at least a projective
line bundle, and so, corresponding to $K_\cS^{1/2}T\cS$, which is trivial on
minitwistor lines and has trivial determinant, there should be a canonical
solution of the $\SL(2,\R)$ Einstein-Weyl Bogomolny equations on any
Einstein-Weyl space, yielding a general formula for the Hitchin-LeBrun
construction on any Einstein-Weyl space. In the final section I shall find
this canonical solution directly.

\section{Einstein-Weyl structures and conformal submersions}\label{EW}

If $g$ is a Riemannian metric and $D=D^g+\omega$ is Einstein-Weyl, then it is
well known that $D^g-\omega$ is Einstein-Weyl if and only if $\omega$ is dual,
with respect to $g$, to a conformal vector field.  If $\omega$ is also
divergence-free with respect to $g$, i.e., $g$ is a \emphdef{Gauduchon gauge}
for $D$, then $\omega$ is dual to a Killing field of $g$. (See
e.g.~\cite{CP1,CP2,Gauduchon4,PS1} for more information on Einstein-Weyl
geometry.)

\begin{prop} Suppose that $(M,\conf)$ is a conformal manifold and
$(D^\sd,D^\asd)$ are compatible Einstein-Weyl structures on $M$. Define a
$1$-form $\theta:=\frac12(D^\sd-D^\asd)$ and a Weyl derivative
$D:=\frac12(D^\sd+D^\asd)$ \textup(the barycentre\textup) so that
$D^\sdasd=D\pm\theta$. Then on the open set where $\theta$ is nonvanishing,
$\xi=\theta/|\theta|$ is a conformal submersion with minimal Weyl derivative
$D^0=2D-D^{|\theta|}=D-|\theta|^{-1}D|\theta|$.

\smallbreak
\noindent\textup{(Here $D^{|\theta|}$ is the exact Weyl derivative defined by
$D^{|\theta|}|\theta|=0$. Since $|\theta|$ is a section of $L^{-1}$, this
means that $D^{|\theta|}=D+|\theta|^{-1}D|\theta|$.)}
\proofof{prop} The standard formula for the dependence of the (normalized)
Ricci tensor on the Weyl derivative gives
\begin{equation*}
r^{D^\sdasd}=r^D\mp D\theta+\theta\tens\theta-\tfrac12|\theta|^2\iden
\end{equation*}
and hence $D\theta=\frac12(r^{D^\asd}-r^{D^\sd})$. Since $D^\sdasd$ are both
Einstein-Weyl, $\sym_0 D\theta = 0$ and so one can write $D\theta = \sigma\,
\iden + F$ where $\sigma$ is a section of $L^{-2}$ and $F$ is a skew
endomorphism of weight $-2$. Direct calculations give
\begin{align*}
D^0 &= D - |\theta|^{-1} ( \sigma \xi-F(\xi) )\\
D^{|\theta|} &= D + |\theta|^{-1} ( \sigma \xi-F(\xi)).
\end{align*}
One now readily checks that $\xi$ is a conformal submersion.
\end{proof}
The following diagram in the affine space of Weyl derivatives summarizes
this Proposition, where $\gamma=|\theta|^{-1}D|\theta|
=|\theta|^{-1} ( \sigma \xi-F(\xi))$.
\begin{diagram}[height=1.6em,width=1.6em]
   &           &D^\sd&           &\\
   &\ruTo  &\uTo^\theta&\rdTo    &\\
D^0&\rTo^\gamma&D&\rTo^{\gamma\;\;}&D^{|\theta|}\\
   &\rdTo  &\uTo^\theta&\ruTo    &\\
   &           &D^\asd&          &
\end{diagram}
In particular $D^0$ is exact iff $D$ is exact. This holds, for instance
if $D^\sdasd$ are both Levi-Civita derivatives of Einstein metrics.

I now specialize to four dimensions and Einstein-Weyl structures with selfdual
Faraday curvature. On a selfdual conformal $4$-manifold, Einstein-Weyl
structures necessarily have selfdual Faraday curvature~\cite{DMJC2} and are
either Einstein or locally hypercomplex~\cite{PS1}.

\addtocounter{thmcounter}{1}
\begin{thm}\label{biEW} Let $D^\sdasd$ be Einstein-Weyl structures
whose Faraday curvatures $F^{D^\sdasd}$ are selfdual. Then, with the
notation of the previous proposition, $\xi$ is a selfdual conformal
submersion, and also $\alt(D^0\tsum D^{sdasd})\xi$ is selfdual.
\proofof{thm} Since $F^{D^\sd}$ and $F^{D^\asd}$ are both selfdual, so is
$F^D$, and hence so is $F^0$, since $D^{|\theta|}$ is exact. This means that
$\xi$ is a selfdual conformal submersion. Furthermore $F$ is selfdual, since
it is equal to $\frac14(F^{D^\sd}-F^{D^\asd})$. Consequently
$\omega=-2|\theta|^{-1}F(\xi)$ and so $D=D^{sd}+|\theta|^{-1}\sigma\xi$.
Therefore $D^\sdasd$ both differ from $D^{sd}$ by a vertical $1$-form, and so
the antiselfdual part of $(D^0\tsum D^\sdasd)\xi$ vanishes.
\end{proof}
In the language of the previous section, the final condition of $\xi$ means
that $\xi$ is Legendrian if $D^\sdasd$ is the Levi-Civita derivative of an
Einstein metric, and triholomorphic if $D^\sdasd$ is the Obata derivative of a
hypercomplex structure.

Another picture in the affine space of Weyl derivatives may be helpful.
\begin{diagram}[width=2em,height=1.5em,tight]
D^0&\rTo^{\frac12\omega}&D^{sd}    &\rTo^{\frac12\omega}&D^B\\
   &\rdTo_\gamma&\dTo_{|\theta|^{-1}\sigma\xi}&              &\\
   &                    &D         &      &\dTo_{2|\theta|^{-1}\sigma\xi}\\
   &            &                  &\rdTo_\gamma             &\\
   &                     &         &                    &D^{|\theta|}
\end{diagram}

When both Einstein-Weyl structures are (locally) hypercomplex $\sigma=0$ and
so $D^B=D^{|\theta|}$ and the Einstein-Weyl structure on $B$ is Einstein.
Indeed, since the Einstein-Weyl quotient $B$ admits two hyperCR structures, it
must be the round $3$-sphere metric~\cite{GT}. Such bi-hypercomplex
structures in four dimensions have been studied already by Apostolov and
Gauduchon [private communication] and they have an elegant construction of the
conformal submersion that arises in this case, which I shall briefly describe.

If the hypercomplex structures corresponding to $D^\sd$ and $D^\asd$ are
$(I^\sd_1,I^\sd_2,I^\sd_3)$ and $(I^\asd_1,I^\asd_2,I^\asd_3)$, then since
both give oriented orthonormal frames for $L^2\Lambda^2_-T\dual M$, they are
related by an $\SO(3)$-valued function: $I^\sd_i=A_{ij}I^\asd_j$.  Applying
$D$ to this equation gives $dA_{ij}(X)I^\asd_j=-2\abrack{\theta\skwend
X,A_{ij}I^\asd_j}$, which implies, after some manipulations, that
$dA_{ij}(X)=-A_{ik}\eps_{jk\ell}\theta(I^\asd_\ell X)$. In particular,
$dA_{ij}(\xi)=0$ so this map to $\SO(3)$ factors through the conformal
submersion $\xi$. To see that the map to $\SO(3)$ actually \emph{is} the
conformal submersion $\xi$, one computes
$dA_{ij}(X)dA_{ij}(X)=2(|\theta|^2|X|^2-\theta(X)^2)$.

Conversely, by Theorem II, any selfdual conformal submersion over the round
$3$-sphere metric is bi-triholomorphic with respect to a bi-hypercomplex
structure. The generalized monopole equations of Proposition~\ref{evolve} in
the gauge $\gmw=|\theta|=\mu_{S^3}^{-1}$ reduce to $\stB\dot A\gmw=dA+\dot
A\wedge A$, which have been obtained independently by Belgun and Moroianu.

\section{Selfdual Einstein metrics and hypercomplex structures}\label{SDE}

In this section I study selfdual conformal $4$-manifolds $M$ admitting a
compatible Einstein metric with nonzero scalar curvature and a compatible
hypercomplex structure.  All such structures are obtained by applying the
Hitchin-LeBrun construction to a hyperCR Einstein-Weyl space $B$ and the
Einstein metric will be found explicitly in terms of the Einstein-Weyl
structure on $B$ and a special solution of the affine monopole equations.

The work of the previous section shows that if $D^g=D-\theta$ is the
Levi-Civita derivative of the Einstein metric and $D^{Ob}=D+\theta$ is the
Obata derivative of the hypercomplex structure, then $\xi=\theta/|\theta|$ is
a Legendrian triholomorphic selfdual conformal submersion. The main goal in
this section is to add one more adjective to this list and prove that the
submersion is affine.

To do this, a section $\gmw_1$ of $L^{-1}$ must be found with
$D^0\gmw_1 = F^0(\xi)$: $\gmw_1$ is then linear part of the affine monopole
$\gmw$ and the affine structure is given by $D^0_\xi+\gmw_1$.

\begin{prop} Write $D^{Ob}=D^{sd}-\kappa\xi$. Then $D^0(2\kappa)=F^0(\xi)$.
\proofof{prop} First recall that $\kappa$ is basic, i.e., $D^0_\xi\kappa=0$.
Now $D^0$ and $D^{Ob}$ are gauge equivalent, since $D^g$ and $D^{|\theta|}$ are
both exact. This implies that
\begin{equation*}
0=d(\omega-2\kappa\xi)(\xi)
=d\omega(\xi)+D^0(2\kappa)=F^B(\xi)-F^0(\xi)+D^0(2\kappa)
\end{equation*}
which proves the proposition since $F^B(\xi)=0$.
\end{proof}
This shows that the linear part of the affine monopole is twice the $\kappa$
monopole of the hyperCR space $B = M/\xi$, and so it satisfies
${\stB D^B(2\kappa)}=F^B$~\cite{GT,CP2}.  Fix a gauge on $B$ so that the
Einstein-Weyl structure is given by a metric $g\low_B$ and a $1$-form
$\omB$. Then $\gmw_1=2\kappa$ and one can take $A_1=\omB$.
In order to fix the affine gauge on $M$ note that $\sigma$ is a nonzero
constant multiple of the scalar curvature of $g$ and so $D^g\sigma=0$ since
$g$ is Einstein. In particular $\sigma$ is nonvanishing by assumption.

\begin{prop} Define $D^\af=D^0+2\kappa\xi$. Then
$D^\af_\xi(\sigma^{-1}|\theta|)=-2$ so $\sigma^{-1}|\theta|$ is an affine
section of $L^1$ with respect to the affine structure.
\proofof{prop} Observe that $(D^0-2\theta)(\sigma^{-1}|\theta|)=0$, i.e.,
$D^{\sigma^{-1}|\theta|}=D^0-2\theta$ and so
$D^\af=D^{\sigma^{-1}|\theta|}+2\theta+2\kappa\xi
=D^{\sigma^{-1}|\theta|}-2|\theta|^{-1}\sigma\xi$. Hence $D^\af(\sigma^{-1}
|\theta|)=-2\xi$.
\end{proof}
One more diagram in the affine space of Weyl derivatives may again clarify the
situation, where the Weyl derivative $\hat D$ with $\hat
D\bigl(|\theta|^{-1})=2\xi$ has been introduced for completeness.
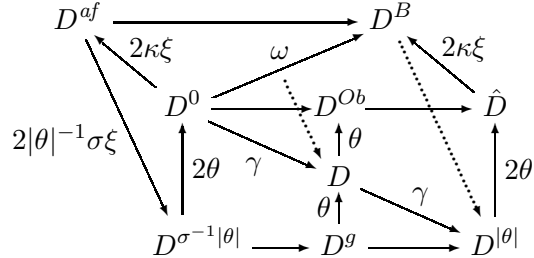
\newdiagramgrid{complicated}{1,1,1,1,1,1,1,1}{1.5,1.5,1.2,1.2,1.2,1.2}
\begin{diagram}[midvshaft,height=1em,width=1.8em,grid=complicated]
D^\af&         &          &    &\rTo&                   &D^B&        &\\
 &\rdTo(2,6)_{\raise4pt\hbox{$2|\theta|^{-1}\sigma\xi\;$}}
  \luTo^{\lower3pt\hbox{$\,2\kappa\xi$}}& &&
 &\ruTo(4,2)^{\lower2pt\hbox{$\omega$}}&
 &\rdDotsto(2,6)\luTo^{\lower3pt\hbox{$\,2\kappa\xi$}}&\\
   &           &D^0       &    &\rTo&\rdDotsto(1,3)[crab=-3pt]\!\!D^{
                          \hbox to6pt{$\scriptstyle Ob$}}& &\rTo&\hat D\\
   &           & &\rdTo(3,2)_\gamma&&   \uTo_\theta      & &    &\\
   &         &\uTo_{2\theta}&  &    &              D     & & &\uTo_{2\theta}\\
   &           &          &    &    &   \uTo^\theta&\rdTo(3,2)^{\gamma\;}& &\\
  & &\quad D^{\sigma^{-1}|\theta|}&&\rTo&          D^g   & &\rTo&D^{|\theta|}
\end{diagram}

It follows from the above proposition that the affine coordinate $t$ can be
taken to be a multiple of $\sigma^{-1}|\theta|$ and so
$|\theta|^{-1}\sigma\xi$ is a multiple of $t^{-1}(dt+t A_1+A_0)$. Now note
that $D^B=D^{|\theta|}-2|\theta|^{-1}\sigma\xi$ and so
$-2d(|\theta|^{-1}\sigma\xi)=F^B$, $A_0=0$ and
$|\theta|^{-1}\sigma\xi=-\frac12(t^{-1}dt+\omB)$. The remaining gauge
freedom can be fixed by taking $\gmw_0=\mu_{\!B}^{-1}$.

Certainly $(1+2t\muB\kappa,t\omB)$ is a solution of the affine monopole
equations on a hyperCR Einstein-Weyl space, yielding a selfdual space with a
hypercomplex structure by Theorem II. The Einstein gauge is the barycentre of
the $|\theta|^{-1}$ gauge and the $\sigma^{-1}|\theta|$ gauge (both of which
have been identified in terms of the affine structure) and one easily checks
that this encodes the Einstein equation, since $D^{Ob}$ is Einstein-Weyl and
$\sym_0 D\theta=0$. Hence, writing $\mu_t=t\muB$, the following theorem is
obtained.

\addtocounter{thmcounter}{1}
\begin{thm}\label{einstein} Let $(\conf\low_B,D^B)$ be
a $3$-dimensional hyperCR Einstein-Weyl space with twist $\kappa$. Then the
metric
\begin{equation*}
g=\frac1{\mu_t^2}\bigl( (1+2\mu_t\kappa)\pi^*\conf\low_B
+(1+2\mu_t\kappa)^{-1} (D^B\mu_t)^2\bigr)
\end{equation*}
on $\pi\colon M\to B$ is selfdual Einstein with nonzero scalar curvature and
admits a compatible hypercomplex structure, where $\mu_t$ is a section of
$L^1$ increasing along the fibres. Any selfdual Einstein metric with a
compatible hypercomplex structure arises in this way.
\end{thm}
\noindent Note that $\conf\low_B$ arises as a conformal infinity at
$\mu_t=0$.

Explicit selfdual Einstein metrics can be found by applying this construction
to explicit hyperCR Einstein-Weyl spaces. To the best of my knowledge, the
known examples are the Einstein-Weyl spaces with geodesic symmetry, and the
hyperCR-Toda spaces. In the former case note that the twist of the hyperCR
structure is \emphdef{minus} the twist of the geodesic symmetry and so the
examples of section~\ref{affine} (with $h=f$) are reobtained~\cite{CT}---these
are the Pedersen metrics~\cite{HP2} when $h$ is constant.

On the other hand, the hyperCR-Toda spaces yield new selfdual Einstein metrics
with compatible hypercomplex structures and no continuous symmetries:
\begin{align*}
g&=\frac1{t^2}
\bigl(H(\,|z+h|^2g\low_{S^2}+dz^2)+H^{-1}(dt+t\omB )^2\bigr)\\
\tag*{where}H&= 1+\frac{i(h-\overline h)t}{|z+h|^2},\qquad
\omB =-\frac{2z+h+\overline h}{|z+h|^2}dz
\end{align*}
and $h$ is holomorphic on an open subset of $S^2$.

HyperCR Einstein-Weyl spaces are not well understood: the results of this
section perhaps provide motivation for further investigations. Alternatively,
one may study hypercomplex selfdual Einstein $4$-manifolds directly. Along
these lines Apostolov and Gauduchon have investigated compatible selfdual
complex structures and obtained a nice characterization of the Pedersen
metrics.

\section{The Hitchin-LeBrun construction}

There is an interesting gauge transformation one can apply to the
metric of Theorem IX: on replacing $\mu_t$ by the new projective coordinate
$\mu_t/(1-\mu_t\kappa)$, the metric becomes:
\begin{equation*}
g= \frac1{\mu_t^2}\Bigl(\bigl(1-\mu_t^2\kappa^2\bigr)\pi^*\conf\low_B
+ \bigl(1-\mu_t^2\kappa^2\bigr)^{-1}
\bigl(D^B\mu_t+\mu_t^2 D^B\kappa)\bigr)\Bigr).
\end{equation*}
Now on a hyperCR Einstein-Weyl space, $D^B\kappa=-\frac12{\stB F^B}$ and
$\kappa^2=\frac16\scal^B$ (see~\cite{GT}), so this form of the metric makes
sense for any Einstein-Weyl space. In this final section, I prove the
following theorem.

\addtocounter{thmcounter}{1}
\begin{thm}\label{final} Let $(\conf\low_B,D^B)$
be an arbitrary $3$-dimensional Einstein-Weyl structure with Faraday curvature
$F^B$ and scalar curvature $\scal^B$. Then
\begin{equation*}\notag
g=\bigl(1-\tfrac16\mu_t^2\scal^B\bigr)\mu_t^{-2}\pi^*\conf\low_B
+\bigl(1-\tfrac16\mu_t^2\scal^B\bigr)^{-1}
\bigl(\mu_t^{-1}D^B\mu_t-\tfrac12\mu_t{\stB  F^B}\bigr)^2
\end{equation*}
is a selfdual Einstein metric of nonzero scalar curvature, with a Legendrian
selfdual conformal submersion over $B$. Here $\mu_t$ is a section of $L^1$
increasing along the fibres of the conformal submersion $\pi\colon M\to B$ and
the conformal structure $\conf\low_B$ is the conformal infinity at
$\mu_t=0$. Any such selfdual Einstein metric arises in this way.
\end{thm}
Strictly speaking, the above metric is only positive definite for
$\mu_t^2\scal^B<1$ and in this region the scalar curvature is negative. For
$\mu_t^2\scal^B>1$ the negation of the above metric is positive definite and
has positive scalar curvature.

Note that the Einstein metric can be written in a gauge $\muB$ by writing
$\mu_t=t\muB$, $g\low_B=\mu_{\!B}^{-2}\conf\low_B$ and
$D^B=D^{\muB}_{\vphantom{T}}+\omB$. Then
\begin{equation*}\notag
g= \frac1{t^2}\Bigl(\bigl(1-\tfrac16t^2\mu_{\!B}^2\scal^B\bigr) g\low_B
+ \bigl(1-\tfrac16t^2\mu_{\!B}^2\scal^B\bigr)^{-1}
\bigl(dt + t \omB  - \tfrac12 t^2\muB {\stB  F^B}\bigr)^2\Bigr).
\end{equation*}

The theorem is proven using the projective monopole equations. The
conformal structure is determined by a canonical $\SL(2,\R)$ monopole
on $L^{1/2}_B\dsum L^{-1/2}_B$. In a gauge $\muB$, the Higgs field and
connection $1$-form are given by:
\begin{equation*}
\gmw=\begin{pmatrix}0& \mu_{\!B}^{-1}\\ \tfrac16\muB \scal^B& 0\end{pmatrix},
\qquad A=\begin{pmatrix}\tfrac12\omB & 0\\
\tfrac12\muB {\stB  F^B}& -\tfrac12\omB \end{pmatrix}.
\end{equation*}
The connection is therefore $D^B+{\stB  F^B}$, with ${\stB  F^B}$
acting from $L^{1/2}_B$ to $L^{-1/2}_B$, while the Higgs field is $1+\scal^B$
in $L^{-1}_B\tens
\bigl(\Hom(L^{-1/2}_B,L^{1/2}_B)\dsum\Hom(L^{1/2}_B,L^{-1/2}_B)\bigr)$.

The Einstein-Weyl Bogomolny equations are:
\begin{align*}
-\tfrac16{\stB (D^B(\muB \scal^B)-\omB \muB \scal^B)}
&=-\tfrac12d(\muB {\stB F^B})+\tfrac12\omB \wedge\muB 
{\stB F^B},\\
\tfrac12 F^B&=\tfrac12d\omB ,\\
D^B(\mu_{\!B}^{-1})+\omB \mu_{\!B}^{-1}&=0.
\end{align*}
Since $D^B\muB =\omB \muB $ and $F^B=d\omB $, the only
nontrivial equation is the first one, which reduces to
$\frac16 D^B\scal^B=-\frac12{\stB d^B{\stB F^B}}=\frac12\delta^BF^B$,
where $\delta^B$ is the twisted exterior divergence. This equation is
satisfied automatically: it is the differential Bianchi identity
for the Weyl connection in Einstein-Weyl geometry~\cite{Gauduchon4,CP1}.

It follows from Theorem VI that the metric $g$ is selfdual with a selfdual
projective conformal submersion. It remains, therefore, to prove that $g$ is
Einstein and $\xi$ is Legendrian. To do this, observe that the Levi-Civita
derivative of $g$ is the barycentre of the Levi-Civita derivatives of $\mu_t$
and $t\gmw^{-1}$. Simple calculations show that these are given by
\begin{align*}
D^B-\mu_t^{-1}D^B\mu_t&=
D^B-(1-\tfrac16\mu_t^2\scal^B)\mu_t^{-1}\xi-\tfrac12\mu_t{\stB F^B}\\
\tag*{and}
D^\gmw-t^{-1}dt&=
D^0-(1+\tfrac16\mu_t^2\scal^B)\mu_t^{-1}\xi+\tfrac12\mu_t{\stB F^B}.
\end{align*}
Hence $D^g=D^{sd}-\mu_t^{-1}\xi$ (which will give the Legendrian property) and
consequently,
\begin{equation*}
r^g=r^{sd}+D^B(\mu_t^{-1})\tens\xi+\mu_t^{-1}(D^0\tsum D^{sd})\xi
+\mu_t^{-2}(\xi\tens\xi-\tfrac12\iden),
\end{equation*}
where I have written $D^{sd}=D^B\tsum D^0\tsum D^{sd}$ on
$T\dual M=L^{-1}\tens L^{-1}\tens TM$. This simplifies to give
\begin{equation*}
r^g=r^{sd}-\tfrac12{\stB F^B}\tens\xi+\tfrac16\scal^B\xi\tens\xi
-\tfrac12\mu_t^{-2}\iden+\mu_t^{-1}(D^0\tsum D^{sd})\xi.
\end{equation*}
Finally substituting from Proposition~\ref{ricci} yields $\sym
r^g=\frac12(\frac16\scal^B-\mu_t^{-2})\iden$, and so $g$ is Einstein, with
scalar curvature $-12\mu_t^{-2}(1-\frac16\scal^B\mu_t^2)$. This completes
the proof of Theorem X.

Examples arising from hyperCR Einstein-Weyl spaces have already been
discussed. One source of further examples are the Ward-Toda
spaces~\cite{Ward,DMJC4}
\begin{align*}
g&=(V_\rho^2+V_\eta^2)(d\rho^2+d\eta^2)+d\psi^2\\
\omega&=\frac{2V_\rho V_\eta\,d\eta+(V_\rho^2-V_\eta^2)d\rho}
{\rho(V_\rho^2+V_\eta^2)}
\end{align*}
where $V$ is an axially symmetric harmonic function: $(\rho V_\rho)_\rho+\rho
V_{\eta\eta}=0$. These spaces admit a symmetry generated by
$\partial/\partial\psi$ and hence so will their Hitchin-LeBrun metrics.
Already in these examples, the Faraday and scalar curvatures are quite
formidable, so these Einstein metrics are not at all simple. Nevertheless,
they can be made completely explicit, and will undoubtedly repay further
study.

\end{document}